\documentclass[twocolumn,english,aps,pre,eqsecnum]{revtex4-2}
\usepackage[T1]{fontenc}
\usepackage[latin9]{inputenc}
\usepackage{geometry}
\usepackage{float}
\usepackage{fancybox}
\usepackage{calc}
\usepackage{bm}
\usepackage{amsmath}
\usepackage{amssymb}
\usepackage{graphicx}

\makeatletter

\providecommand{\tabularnewline}{\\}

\makeatother

\usepackage{babel}
\begin{document}
\title{A midpoint projection algorithm for stochastic differential equations
on manifolds}
\author{Ria Rushin Joseph$^{1}$, Jesse van Rhijn$^{1,2}$, Peter D. Drummond$^{1}$}
\affiliation{$^{1}$Centre for Quantum Science and Technology Theory,~\\
Swinburne University of Technology, Melbourne, Victoria, Australia}
\affiliation{$^{2}$University of Twente, Enschede, The Netherlands}
\begin{abstract}
Stochastic differential equations projected onto manifolds occur in
physics, chemistry, biology, engineering, nanotechnology and optimization,
with interdisciplinary applications. Intrinsic coordinate stochastic
equations on the manifold are often computationally impractical, and
numerical projections are useful in many cases. We show that the Stratonovich
interpretation of the stochastic calculus is obtained using adiabatic
elimination with a constraint potential. We derive intrinsic stochastic
equations for spheroidal and hyperboloidal surfaces for comparison
purposes, and review some earlier projection algorithms. In this paper,
a combined midpoint projection algorithm is proposed that uses a midpoint
projection onto a tangent space, combined with a subsequent normal
projection to satisfy the constraints. Numerical examples are given
for a range of manifolds, including circular, spheroidal, hyperboloidal,
and catenoidal cases, as well as higher-order polynomial constraints
and a ten-dimensional hypersphere. We show that in all cases the combined
midpoint method has greatly reduced errors compared to methods using
a combined Euler projection approach or purely tangential projection.
Our technique can handle multiple constraints. This allows manifolds
that embody several conserved quantities. The algorithm is accurate,
simple and efficient. An order of magnitude error reduction in diffusion
distance is typically found compared to the other methods, with reductions
of several orders of magnitude in constraint errors.
\end{abstract}
\maketitle

\section{Introduction}

Diffusion processes and stochastic trajectories on curved manifolds
inside Cartesian spaces have many applications. The diffusion of atoms
or molecules on many types of curved surface occurs in surface physics,
biophysics, catalysis, biochemistry, cell biology, and nanotechnology
\citep{cherry1979rotational,Brillinger1997,lin2004dynamics,sbalzarini2006simulations,novak2007diffusion,gusak2010diffusion,klaus2016analysis,adler2019conventional,castro2018active}.
Curved-space diffusion also arises in more abstract problems in general
relativity \citep{smerlak2012diffusion}, imaginary-time path-integrals
\citep{kosztin1996introduction} and in quantum field theory \citep{mattis1998uses}. 

Other related problems include systems of stochastic differential
equations (SDEs) with conserved quantities \citep{zhou2016projection,albeverio1995remark}.
These are used for sampling purposes in studying classical free energies
\citep{ciccotti2008projection}, and molecular dynamics. Multidimensional
constraints are also found in engineering systems \citep{schonlau1998global,lee2006mining},
such as in robotics \citep{meng2017iterative,carius2021constrained,stilman2010global},
where joints may have restricted mobility due to their mechanical
design, and experience stochastic forces. 

In this paper a hybrid midpoint projection algorithm is introduced
for general projected SDEs which combines both tangential and normal
projections. A related approach uses an Euler method for each tangential
step \citep{holyst1999diffusion}. Here we generalize this approach
to a more accurate midpoint projection algorithm for curved surface
stochastic equations.

The resulting algorithm is robust and easily implemented, treating
arbitrary drift, diffusion, constraint type and projected manifold
dimension. Numerical examples are given, demonstrating excellent error
performance both in measures of diffusion distance and in satisfaction
of the constraints, which can give rise to large global errors with
other methods.

Phase-space stochastic processes \citep{chierchia1994drift,ngai1991toward,thomson1972turbulent}
as well as quantum phase-space methods, for example representations
\citep{Ria:2018_Majorana,Riashock2018} of Majorana physics \citep{beenakker2013search,wilczek2009majorana}
also lead to curved space Fokker-Planck \citep{graham1977lagrangian,Graham1977Covariant,RiaJoseph2021interacting}
and stochastic equations \citep{Gardiner1997}. These are some of
the many applications of embedded stochastic processes. 

Constraint projection occurs in molecular dynamics calculations. For
these, specialized techniques are known. Such earlier algorithms are
often restricted to these specific problems \citep{andersen1983rattle},
or are optimized for steady-state distribution sampling, which may
not give optimal performance for dynamics \citep{Lelievre2010free,laurent2021order}. 

Due to the diversity of the applications, here we develop a general
technique applicable to any class of projected stochastic differential
equation, including multiplicative noise and arbitrary projections.

Projected diffusion processes can sometimes be treated using stochastic
equations for a lower dimensional set of coordinates on the manifold,
for example by using spherical polar coordinates on a spherical surface.
These can be termed implicit projections: the manifold is implicit
in the definition of the coordinates. The embedding theorems \citep{whitney1936differentiable,whitney1944,nash1956imbedding}
of Riemannian geometry show that surfaces of arbitrary smooth curvature
are equivalent to an appropriate projection. The resulting stochastic
equations usually involve multiplicative noise, because the noise
coefficients depend on the metric tensor of the manifold. Several
numerical algorithms are known for treating such multiplicative stochastic
differential equations (SDE)s \citep{Drummond1990,Kloeden1992,milstein1994numerical}.

There are problems with implicit projections \citep{bertalmio2001variational},
especially for closed surfaces. The projected equations may have singularities
in the noise and drift terms, due to the analyticity properties of
maps \citep{hawley1950theorem} and the Poincare-Hopf theorem \citep{poincare1885courbes,hopf1927vektorfelder}.
Stochastic numerical algorithms may not converge near singularities.
Also, the definition of an intrinsic coordinate system \citep{emery1990two,calhoun2010finite}
on curved surfaces can be a complex problem. There is generally no
simple reduced coordinate system for a biological cell wall or other
low-dimensional structures in biology and engineering, due to their
complex geometry.

The approach used here is to define a stochastic process in the original
coordinates, whose projection satisfies the restrictions of an embedded
manifold. This overcomes singularity and variable change issues. Even
if the stochastic equation is known on the manifold, it can be useful
to employ projections in order to eliminate singularities. The drawback
is that discretization can lead to errors from finite step-sizes which
may take the stochastic trajectory off the manifold. Such errors may
grow rapidly in time as the trajectory moves progressively away from
the manifold. These additional global errors are not present in the
algorithm proposed here, since the constraints are satisfied with
a normal projection after each step.

We first review the theory of manifolds embedded in a Euclidean space.
An explanation of how projections can occur physically is obtained
through adiabatic elimination with constraint potentials. This is
then applied to projections of stochastic processes. The continuous
time limit of a constrained stochastic equation is derived as a Stratonovich
equation, which is shown to be valid under conditions of finite bandwidth
noise and fast projection timescales. Algorithms for multi-dimensional
tangential and normal projections are derived. These methods are used
to generalize some previous proposals for solving projected stochastic
equations. 

In order to compare the methods, numerical implementations are carried
out. The first examples treat diffusion on circular surfaces, as well
as spheroidal and hyperboloidal surfaces in a three-dimensional Cartesian
space. We use exact solutions or a semi-implicit intrinsic method
using very small time-steps to generate reference results where possible.
Simulations using hybrid Euler \citep{holyst1999diffusion} and tangential
methods \citep{armstrong2019optimal} are also carried out, with errors
calculated for measures of distance diffused, and for conservation
of constraints. We show that the hybrid midpoint projection algorithm
has lower errors in both categories than the other methods. 

An outline of the paper is as follows. In Section \ref{sec:Manifolds,-constraints-and},
we review the geometric background of manifolds, constraints and projections.
In Section \ref{sec:Adiabatic-elimination-for} adiabatic elimination
is used to justify the interpretation of the projected equation as
a Stratonovich equation. Intrinsic stochastic equations are obtained
for general spheroidal and one-sheet hyperboloidal surfaces. Section
\ref{sec:Stochastic-projection-algorithms} gives a discussion of
numerical algorithms and the different types of projections. Section
\ref{sec:Numerical-examples} gives numerical examples and error comparisons
with either exact or intrinsic diffusion results for circular, catenoidal,
spheroidal and hyperboloidal surfaces embedded in Euclidean space.
Higher-order polynomial surfaces and hyperspheres are also investigated
to model non-quadratic and higher dimensional surfaces. Section \ref{sec:Summary}
summarizes the results. In the Appendix we give the intrinsic SDE
variables and derivations for spheroidal and hyperboloidal surfaces.

\section{Manifolds, constraints and projections\label{sec:Manifolds,-constraints-and}}

This paper treats numerical algorithms for projected stochastic equations.
We utilize Stratonovich equations \citep{stratonovich1966new,Gardiner1997},
which occur in physics as the broad-band limit of a finite bandwidth
random process. These follow standard variable change rules. Projections
using Ito stochastic calculus \citep{armstrong2019optimal}, are not
treated here, and have different properties. While this section reviews
standard results, they are useful to define notation.

An unprojected Stratonovich SDE on a Euclidean space $\mathbb{R}^{n}$
defines a real vector trajectory $\bm{x}=\left[x^{1},\ldots x^{n}\right],$
where
\begin{equation}
\dot{\bm{x}}=\bm{a}\left(\bm{x},t\right)+\bm{B}\left(\bm{x},t\right)\bm{\xi}\left(t\right).\label{eq:Original-Strat}
\end{equation}

The standard problem is to solve this over a finite interval $\left[0,T\right]$,
giving a time dependent functional probability density of trajectories
$\mathcal{P}\left(\left[\bm{x}\left(t\right)\right]\right)$ conditioned
on a specified initial probability distribution $P\left(\bm{x},t=0\right)$.
It is usual to focus on the marginal probability densities $P\left(\bm{x},t\right)$.

Here $\bm{a}\left(\bm{x},t\right)$ is a real vector function of $\bm{x}$
of dimension $n$, $\bm{B}\left(\bm{x},t\right)$ is a real $n\times s$
dimensional matrix function, where we omit the arguments when there
is no ambiguity. The real $s$-dimensional Gaussian noise vector $\bm{\xi}\left(t\right)$
physically has a finite correlation time $T_{c}$, and is delta-correlated
in the broad-band limit of small $T_{c}$, so that: 
\begin{equation}
\lim_{T_{c}\rightarrow0}\left\langle \xi^{\sigma}\left(t\right)\xi^{\kappa}\left(t'\right)\right\rangle =\delta^{\sigma\kappa}\delta\left(t-t'\right).\label{eq:noise-correl}
\end{equation}

For a small time interval $\Delta t=t_{1}-t_{0}\gg T_{c}$, we define
a stochastic integral $\Delta w_{i}$ as:
\begin{equation}
\Delta w^{\sigma}=\int_{t_{0}}^{t_{1}}\xi^{\sigma}\left(t\right)dt,
\end{equation}
which has the property that 
\begin{equation}
\left\langle \Delta w^{\sigma}\Delta w^{\kappa}\right\rangle =\delta^{\sigma\kappa}\Delta t.\label{eq:Finite-step-noise-correlation}
\end{equation}

A single discrete step in the stochastic trajectory is $\Delta\bm{x}=\bm{x}\left(t_{1}\right)-\bm{x}\left(t_{0}\right).$
In the Stratonovich interpretation of stochastic calculus, the midpoint
$\bar{\bm{x}}$ is the location where the derivatives are evaluated,
where
\begin{equation}
\bar{\bm{x}}=\left(\bm{x}\left(t_{1}\right)+\bm{x}\left(t_{0}\right)\right)/2.\label{eq:midpoint}
\end{equation}
Introducing the corresponding mean time $\bar{t}=\left(t_{1}+t_{0}\right)/2$,
the broad-band limit of Eq (\ref{eq:Original-Strat}) is defined \citep{stratonovich1966new}
as the limiting behavior of the following implicit stochastic difference
equation
\begin{equation}
\Delta\bm{x}=\bm{a}\left(\bar{\bm{x}},\bar{t}\right)\Delta t+\bm{B}\left(\bar{\bm{x}},\bar{t}\right)\Delta\bm{w}.\label{eq:Finite-step-SDE}
\end{equation}

In many applications one must project this equation onto a manifold.
The broad-band limit is only taken here after adiabatic elimination
\citep{gardiner1984adiabatic}, and we assume that the time-scales
of the constraint process are much faster than the noise correlation
time, $T_{c}$. Although we only treat real spaces here, these results
are applicable to complex spaces, with minor changes. In this section
we analyze how such manifolds are obtained, and the corresponding
geometrical definitions. 

Since we wish to derive algorithms numerically applicable to this
problem, we use a coordinate-dependent notation rather than a coordinate-free
approach. The results will be expressed in the coordinates of the
enclosing (Euclidean) space. In the examples, where possible, explicit
comparisons are made to intrinsic coordinate methods. 

\subsection{Manifolds}

In the section we review the notation that defines the manifold and
its properties \citep{frankel2011geometry,schuster2003quantum}. We
equip $\mathbb{R}^{n}$ with the Euclidean metric. We use standard
contravariant notation for coordinates $\bm{x}=\{x^{i}\}$, and covariant
notation for derivatives and normals $\bm{n}=\{n_{i}\}$ to allow
the results to be extended to other metrics. Repeated induces will
be summed over in the Einstein convention. In the Euclidean case $x^{i}=x_{i}$,
so  the index placement can be ignored.

An $m$-dimensional sub-manifold is the set of points $\bm{y}\in\mathcal{M}\subset\mathbb{R}^{n}$
for which
\begin{equation}
\bm{f}\left(\bm{y}\right)=0,\label{eq:general-constraint}
\end{equation}
where $\bm{f}$ is a $p=n-m$ dimensional vector constraint function
that defines the projection. We assume that the constraint equations
are linearly independent. Here $\bm{f}=\left(f^{1},\ldots f^{p}\right)$
denotes the $p$ constraints.

Local curvature properties are crucial in determining the properties
of projected stochastic equations. As a result, it is useful to focus
on the classical quadratic geometries. These cover a wide range of
Riemann curvatures and topologies. They include positive and negative
curvature, and different manifold dimensions and co-dimensions. Locally,
most realistic geometries will be comparable to at least one of these
cases, and we treat more general cases both in the analytic theory
and in the examples.

In the simplest examples below, the constraints are quadratic functions,
where $f_{0}^{i}$ is a real constant, $\bm{h}^{i}$ is a real vector,
$\bm{G}^{i}$ is an $N\times N$ real matrix and:
\begin{equation}
f^{i}\left(\bm{y}\right)=f_{0}^{i}+\bm{h}^{i}\cdot\bm{y}+\bm{y}^{T}\bm{G}^{i}\bm{y}.\label{eq:constraint}
\end{equation}

This covers circles, parabolas, spheres, hyper-spheres, spheroids,
and hyperboloids, as well as intersections of these surfaces when
there are multiple constraints. We use these for examples, and choose
cases such that $f^{0}=-1$ and $\bm{h}=0$ for simplicity. The general
algorithms and examples discussed are not limited to these quadratic
constraints, but quadratic constraints are useful because these manifolds
have well-known intrinsic coordinates and diffusion equations for
comparison purposes. 

In some cases it is possible to define $\bm{\phi}\in\mathbb{R}^{m}$
as a new intrinsic coordinate system for the set of points restricted
to the manifold. For example, on the surface of a sphere one may use
polar coordinates. In such cases, a vector mapping function $\bm{\Phi}$
exists where $\bm{\phi}=\bm{\Phi}\left(\bm{y}\right)$ for every Euclidean
coordinate $\bm{y}$ in the original manifold $\mathcal{M}$. 

Even if they exist, such mapping functions lead to singularities for
a compact manifold, and may not be able to be expressed in a closed
form. As a result, we use algorithms that can be expressed in the
original Euclidean coordinates. Intrinsic coordinates are used for
numerical examples, although such comparisons are only available in
special cases where analytic transforms exist.

\subsection{Tangent and normal spaces}

In many cases it is not convenient or practical to use intrinsic coordinates,
and it is necessary to integrate in the original Euclidean frame of
reference. This requires a knowledge of the tangential and normal
spaces. Near to a point $\bm{y}$ in the manifold, from a Taylor expansion
of the $j$-th constraint equation,
\begin{equation}
f^{j}\left(\bm{x}\right)=\sum_{i=1}^{n}\Delta^{i}f_{,i}^{j}\left(\bm{y}\right)+O\left(\left|\bm{\Delta}\right|^{2}\right),
\end{equation}
where $\Delta^{i}=x^{i}-y^{i}$ , $\partial_{i}\equiv\partial/\partial x^{i},$
and the constraint derivative is
\begin{equation}
f_{,i}^{j}\equiv\partial_{i}f^{j}.
\end{equation}

From this condition, $v_{\perp i}^{j}=\partial_{i}f^{j}$ is called
the $j$-th normal vector, $\bm{v}_{\perp}^{j}$, where the lower
index indicates a covariant vector, and the upper index labels the
constraint, $f^{j}$, so that $j=1,\ldots p$. Any orthogonal vector
$\bm{v}^{\parallel}$ such that $\bm{v}^{\parallel}\cdot\bm{v}_{\perp}^{j}=0$
locally satisfies the $j-$th constraint, although it may not satisfy
the other constraints. In general the normal vectors $\boldsymbol{v}_{\perp}^{j}$
may not be orthogonal to each other. They are assumed to be linearly
independent at all locations in order to generate a differentiable
manifold locally isomorphic to $\mathbb{R}^{m}$ without a singular
point \citep{frankel2011geometry}. Singular points may occur with
some choices of multiple projections, but these are unphysical and
can lead to convergence issues.

One can use the gradient vectors of the constraints to define a set
of orthonormal vectors $\boldsymbol{n}^{j}$ that span the $p$-dimensional
space of vectors normal to each constraint. This requires orthogonalization
of the set $\boldsymbol{v}_{\perp}^{j}$, either by Gram-Schmidt or
other methods \citep{trefethen1997numerical}. As a result, the orthonormal
vectors $\bm{n}^{j}$ have the property that, for $i,j=1,\ldots p$:
\begin{align}
\bm{n}^{j}\cdot\boldsymbol{v}_{\perp}^{j} & \neq0\nonumber \\
\bm{n}^{i}\cdot\boldsymbol{n}^{j} & =\delta^{ij}.\label{eq:orthogonal-set-condition}
\end{align}

Similarly, there is a set of orthonormal tangent vectors $\bm{m}_{i}$,
where $i=1,\ldots n$. These define the tangent space of the manifold,
and give an orthonormal basis of $m$ contravariant vectors orthogonal
to $all$ of the normal vectors $\boldsymbol{v}_{\perp}^{j}$, so
that, for $j=1,\ldots p$ and $i=1,\ldots m$:
\begin{align}
\bm{m}_{i}\cdot\boldsymbol{n}^{j} & =\bm{m}_{i}\cdot\boldsymbol{v}_{\perp}^{j}=0\nonumber \\
\bm{m}_{i}\cdot\boldsymbol{m}_{j} & =\delta_{ij}.\label{eq:tangent vectors}
\end{align}

One can define a vector space consisting of all vectors tangent to
$\mathcal{M}$ at any point $\bm{y}\in\mathcal{M}$. This is the tangent
space, $T_{y}\mathcal{M}$, and it is spanned by the basis $\bm{m}_{j}$
for $j=1,\dots m$. The orthogonal complement to the tangent space
is the normal space, spanned by the basis $\bm{n}^{j}$ for $j=1,\dots p$. 

For every vector $\boldsymbol{v}\in\mathbb{R}^{N}$ and manifold coordinate
$\bm{y}\in\mathbb{R}^{m},$ a vector decomposition can therefore be
written as follows \citep{chen2021study}:
\begin{align}
\boldsymbol{v}=\boldsymbol{v}_{\parallel}+\boldsymbol{v}_{\perp},\ \boldsymbol{v}_{\parallel}\in T_{x}\mathcal{M},\ \boldsymbol{v}_{\perp}\in T_{x}^{\perp}\mathcal{M},\label{eq:ss-2}
\end{align}
provided $\boldsymbol{v}_{\parallel}$ , $\boldsymbol{v}_{\perp}$
are the tangential and normal projection of $\boldsymbol{v}\in\mathbb{R}^{m}$
respectively and $T_{y}^{\perp}\mathcal{M}$ is the orthogonal complement
of $T_{y}\mathcal{M}$. If we consider the sphere, the vectors $\boldsymbol{v}_{\parallel}$
perpendicular to the radii are tangents at $\bm{y}$ with a basis
$\bm{m}_{i}$, while the radial vectors $\boldsymbol{v}_{\perp}$
are the normals, with a basis $\boldsymbol{n}^{j}$. 

This can be expressed in terms of projection operators $\mathcal{P}^{\perp}$
and $\mathcal{P}_{\bm{y}}^{\parallel}$ . These satisfy the fundamental
property of mathematical projections on sets that, once projected,
a coordinate doesn't change under further projections, so $\mathcal{P}\left(\mathcal{P}\left(\bm{x}\right)\right)=\mathcal{P}\left(\bm{x}\right)$. 

\subsection{Tangential projections}

Eq. (\ref{eq:tangent vectors}) guarantees that any tangent vector
of form $\bm{y}+\epsilon\bm{m}_{i}$ is in the tangent space. In a
neighborhood of a point $\bm{y}\in\mathcal{M}$ , all constraints
are satisfied up to terms of order $\epsilon^{2}$ since for tangent
vectors $\bm{m}_{i}$ and a small coefficient $\epsilon$:
\begin{equation}
f^{j}\left(\bm{y}+\epsilon\bm{m}_{i}\right)=\epsilon\bm{m}_{i}\cdot\boldsymbol{v}_{\perp}^{j}+O\left(\epsilon^{2}\right)\approx0.
\end{equation}

A general tangential projection of any time-evolution equation is
obtained by projecting an arbitrary derivative vector $\bm{\delta}$
at each location onto the tangent space of the manifold. This type
of projection depends on the location, so we use the notation of $\bm{\delta}^{\parallel}=\Pi_{\bm{y}}^{\parallel}\left(\bm{\delta}\right)$,
for a tangential projection of a derivative vector $\bm{\bm{\delta}}$
at a point on or near the manifold. 

One can obtain the tangential projection $\bm{\delta}^{\parallel}$
in a simple case of a one-dimensional constraint, by removing the
normal component from the vector $\bm{\bm{\delta}}$ using Eq. (\ref{eq:ss-2}),
so that:
\begin{align}
\bm{\delta}^{\parallel} & =\bm{\delta}-\left(\bm{\delta}.\boldsymbol{v}_{\perp}\right)\boldsymbol{v}_{\perp}/\left|\boldsymbol{v}_{\perp}\right|^{2},\nonumber \\
 & \equiv\left(I-\boldsymbol{n}\otimes\boldsymbol{n}\right)\bm{\delta}.
\end{align}
Here $\boldsymbol{n}\otimes\boldsymbol{n}$ is an outer vector product
defined as the matrix $\left[\boldsymbol{n}\otimes\boldsymbol{n}\right]_{ij}=n_{i}n_{j}$,
$\boldsymbol{n}$ is the normal vector at $\bm{y},$ $\boldsymbol{I}$
is the identity matrix, and $\left(\boldsymbol{I}-\boldsymbol{n}\otimes\boldsymbol{n}\right)$
is called the tangential operator \citep{pozrikidis1997introduction}. 

More generally, for multiple constraints, the projection must remove
multiple normal components $\boldsymbol{n}^{j}$, so one must define:
\begin{align}
\bm{\delta}^{\parallel} & =\mathcal{P}_{\bm{y}}^{\parallel}\left(\bm{\delta}\right)\\
 & =\left(I-\sum_{j=1}^{p}\boldsymbol{n}^{j}\otimes\boldsymbol{n}^{j}\right)\bm{\delta}.\nonumber 
\end{align}

\subsection{Normal projections}

In a normal projection, points in a neighborhood of $\mathcal{M}$
are projected onto the nearest point on $\mathcal{M}$. The projection
$\bm{y}=\mathcal{P}^{\perp}\left(\bm{x}\right)$ therefore takes an
initial $\bm{x}\in\mathbb{R}^{n}$, and maps it onto the nearest manifold
coordinate $\bm{y}$. 

For example, in the spherical case, a path in Euclidean coordinates
$\bm{x}$ in $n$-space is projected onto the $n-1$ dimensional hyper-spherical
surface of a ball of radius $r$ around the origin, at every step
on the path. An exact projection function in this case is:
\begin{equation}
\bm{y}=\mathcal{P}^{\perp}\left(\bm{x}\right)=\frac{r\bm{x}}{\left|\bm{x}\right|},
\end{equation}
 which is defined for all $\bm{x}\neq0$. Exact normal projections
are only obtainable in special cases. Hence, we now consider approximate
normal projections valid near a manifold point. 

To obtain the general approach which solves this, we introduce a Lagrange
multiplier vector $\bm{c}$. Given an initial point $\bm{x}$ and
a projected point $\bm{y}$ on the manifold, one must find Lagrange
multipliers $c_{j}\left(\bm{x},\bm{y}\right)$ such that: 
\begin{align}
\bm{y} & =\bm{x}+\sum_{j=1}^{p}c_{j}\left(\bm{x},\bm{y}\right)\boldsymbol{n}^{j}\left(\bm{y}\right),\nonumber \\
f^{j}(\bm{y}) & =0.\label{eq:Projection-2}
\end{align}
Solving for $\bm{y}$ corresponds to minimizing the functional $k(\bm{y})=\left|\bm{y}-\bm{x}\right|^{2}$,
subject to the constraint that $f^{i}\left(\bm{y}\right)=0$. This
is generally a hard computational problem, but can be solved approximately
with a Taylor expansion, valid if $\left|\bm{y}-\bm{x}\right|$ is
small. Since $f^{j}\left(\bm{y}\right)=0$, therefore to first order
in a Taylor series about $\bm{y}$, one finds that
\begin{equation}
f^{j}\left(\bm{x}\right)=\left(\bm{x}-\bm{y}\right)\cdot\boldsymbol{v}_{\perp}^{j}(\bm{y})+O\left(\bm{x}-\bm{y}\right)^{2}.
\end{equation}
This follows since the $\boldsymbol{v}_{\perp}^{j}$ vectors are all
derivatives. From Eq (\ref{eq:Projection-2}), 
\begin{equation}
\bm{y}-\bm{x}=\sum_{j=1}^{p}c_{j}\left(\bm{x},\bm{y}\right)\boldsymbol{n}^{j}\left(\bm{y}\right).
\end{equation}

One must therefore solve the following nonlinear implicit equation:
\begin{equation}
f^{i}\left(\bm{x}\right)=-\sum_{j=1}^{p}M^{ij}\left(\bm{y}\right)c_{j}\left(\bm{x},\bm{y}\right)+O\left(\bm{x}-\bm{y}\right)^{2}.
\end{equation}
In this equation, the normal projection matrix $\bm{M}$ is defined
so that its components are the the inner products of $\boldsymbol{v}_{\perp}^{i}\left(\bm{y}\right)$
with the normal vectors $\boldsymbol{n}^{j}\left(\bm{y}\right)$,
i.e, 
\begin{equation}
M^{ij}\left(\bm{y}\right)=\boldsymbol{v}_{\perp}^{i}\left(\bm{y}\right)\cdot\boldsymbol{n}^{j}\left(\bm{y}\right).\label{eq:normsl-projection-matrix}
\end{equation}
Solving for $c_{j}$, one obtains an implicit equation:
\begin{equation}
c_{i}\left(\bm{x},\bm{y}\right)=-\sum_{j=1}^{p}\left[\bm{M}\left(\bm{y}\right)\right]_{ij}^{-1}f^{j}\left(\bm{x}\right).\label{eq:Implicit-projection}
\end{equation}

Given an initial estimate $\bm{x},$ a normal projection is obtained
by solving Eq. (\ref{eq:Implicit-projection}) for $\bm{c}$, and
hence obtaining $\bm{y}$ using iteration, where the matrix $\bm{M}\left(\bm{x}\right)$
is used instead of $\bm{M}\left(\bm{y}\right)$ to lowest order. In
the algorithms described here, a single iterative step is used because
the trajectories are already close to the manifold, so that:
\begin{equation}
\bm{y}=\mathcal{P}^{\perp}\left(\bm{x}\right)\approx\bm{x}-\sum_{i,j=1}^{p}\boldsymbol{n}^{i}\left(\bm{x}\right)\left[\bm{M}\left(\bm{x}\right)\right]_{ij}^{-1}f^{j}\left(\bm{x}\right).
\end{equation}

In the case of a one-dimensional constraint, since $\left|\boldsymbol{n}^{\perp}\left(\bm{x}\right)\right|=1$,
it follows that $\boldsymbol{v}^{\perp}\left(\bm{x}\right)\cdot\boldsymbol{n}\left(\bm{x}\right)=\left|\boldsymbol{v}^{\perp}\left(\bm{x}\right)\right|$.
For small $\left|\bm{y}-\bm{x}\right|$, one therefore obtains:

\begin{equation}
\bm{y}\approx\bm{x}-\boldsymbol{v}^{\perp}\left(\bm{y}\right)\left|\boldsymbol{v}^{\perp}\left(\bm{x}\right)\right|^{-2}f\left(\bm{x}\right).
\end{equation}

\section{Adiabatic elimination for projected equations \label{sec:Adiabatic-elimination-for}}

To understand the origin of projections of a stochastic equation,
one must ask: what physical or mathematical process restricts the
path to a manifold? Here, we treat a constraining potential, as a
common situation found in physics, and also used for quantum applications
where there are similar issues \citep{schuster2003quantum}. While
this is not universally applicable for all projected stochastic equations,
it explains our choice of stochastic calculus. 

We introduce a local mapping function $\bm{\Phi}$ inside a local
patch where $\bm{\phi}^{\parallel}=\bm{\Phi}\left(\bm{y}\right)$
for Euclidean coordinates $\bm{y}$ in the original manifold $\mathcal{M}$.
For purposes of adiabatic elimination, we extended this locally to
an invertible mapping $\bm{x}=\bm{x}\left(\bm{\phi}\right)$ including
nearby off-manifold points, with a complete set of coordinates $\bm{\phi=}\left[\bm{\phi}^{\parallel},\bm{\phi}^{\perp}\right]$
such that $f^{j}\left(\bm{x}\left(\bm{\phi}\right)\right)=0$ if $\bm{\phi}^{\perp}=0$. 

Specifically, we define $\bm{\phi}^{\perp}$ to correspond to points
in the normal space, where if $\bm{y}$ is a manifold point, then
for nearby points $\bm{x}$ in Euclidean space, $\bm{\phi}^{\perp}$
is the distance along a normal vector:
\begin{equation}
\phi^{\perp,i}\left(\bm{x}\right)=\left(\bm{x}-\bm{y}\right)\cdot\bm{n}^{i}{}^{\perp}\left(\bm{y}\right).\label{eq:normal-variables}
\end{equation}

To return to the spherical example, the coordinate $\phi^{\perp}$
for overall dimension $n=3$ is the radial coordinate in that case.

\subsection{Constraining potentials}

Constraining potentials are essential for constraining particles on
a manifold, and provide a model for the origin of a projected stochastic
equation. We assume that the constraint conditions originate in a
scalar constraint potential 
\begin{equation}
u\left(\bm{x}\right)=\frac{\lambda}{2}\sum_{j=1}^{p}\left[f^{j}\left(\bm{x}\right)\right]^{2}.
\end{equation}
The projected stochastic equation is the adiabatic limit for $\lambda\rightarrow\infty$,
of a continuous stochastic process in which the initial SDE in Eq.
(\ref{eq:Original-Strat}) includes the constraint as a potential,
so that:
\begin{equation}
\dot{x}^{i}=a^{i}\left(\bm{x},t\right)-\partial^{i}u\left(\bm{x}\right)+\sum_{k}B_{\sigma}^{i}\left(\bm{x},t\right)\xi^{\sigma}.\label{eq:constrained SDE}
\end{equation}

Here $\partial^{i}u\equiv\left[\nabla u\right]^{i}$ is the contravariant
derivative, which equals the covariant derivative $\partial_{i}u$
for our metric choice. The drift for motion including the confining
potential is then $\bm{a}_{\lambda}(\bm{x})=\bm{a}-\nabla u\left(\bm{x}\right)$.
The gradient of the potential is,
\begin{equation}
\nabla u\left(\bm{x}\right)=\lambda\sum_{j}\boldsymbol{v}_{\perp}^{j}\left(\bm{y}\right)f^{j}\left(\bm{x}\right).\label{eq:gradient-potential}
\end{equation}
In terms of the local coordinates $\bm{\phi}$, the constrained equation
is:

\begin{equation}
\dot{\bm{\phi}}=\bm{J}\left[\bm{a}_{\lambda}+\bm{B}\bm{\xi}\right].\label{eq:constrained SDE-phi}
\end{equation}
where the Jacobian $\bm{J}$ is defined as:
\begin{equation}
J_{j}^{i}=\frac{\partial\phi^{i}}{\partial x^{j}}.
\end{equation}

From the definition of the normal coordinates in Eq (\ref{eq:normal-variables}),
this implies that for $i>m$, where $m$ is the manifold dimension,
and $\bm{x}\approx\bm{y}$, the Jacobian is equal to the corresponding
normal vector: 
\begin{equation}
J_{j}^{i}=n_{j}^{i}\left(\bm{y}\right).
\end{equation}

The reduction of the number of independent degrees of freedom that
result from a minimization of the potential $u$ is equivalent to
adiabatic elimination of the fast transverse variables in a stochastic
equation, treated next.

\subsection{Adiabatic elimination}

Adiabatic elimination will be used in order to analyze how projections
are applied to stochastic equations. As has been recently pointed
out \citep{armstrong2019optimal}, there can be ambiguity in defining
the which type of stochastic equation is obtained after projection.
This depends on how the projection occurs. We assume that the original
SDE in Eq. (\ref{eq:Original-Strat}) has noise with a finite bandwidth.
Hence, it follows the Stratonovich calculus in the broad-band limit
\citep{stratonovich1966new}. 

The projection is assumed to occur through an adiabatic elimination
process, due to the constraint potential. The resulting diffusion
or noise is generally state-dependent, even if not originally. To
analyze this, we follow a similar method to the work of Gardiner \citep{gardiner1984adiabatic}.
In general, a direct approach to adiabatic elimination is obtained
through dividing up the variables into a ``fast'' and ``slow''
set. In the present context, these are locally in the direction of
the normal and tangential coordinates.

We take a given point $\bar{\bm{\phi}}$ on the manifold such that
$f\left(\bm{x}\left(\bar{\bm{\phi}}\right)\right)=0$. We now consider
a trajectory $\bm{\phi}\left(t\right)$ near $\bar{\bm{\phi}}$, governed
by the constrained stochastic equation in Eq. (\ref{eq:constrained SDE-phi}),
so that 
\begin{equation}
\bm{\Delta}\left(t\right)\equiv\bm{\phi}\left(t\right)-\bar{\bm{\phi}}.
\end{equation}
Changes in the coordinate, $\bm{\Delta}$, are divided into ``slow''
coordinates $\Delta^{\parallel}$ for motion inside the manifold,
and ``fast'' or normal coordinates $\Delta^{\perp}$ for motion
outside the manifold. For a small displacement we obtain:
\begin{equation}
\bm{\Delta}\left(t\right)=\sum_{i=1}^{m}\Delta_{i}^{\parallel}\left(t\right)+\sum_{j=1}^{p}\Delta_{j}^{\perp}\left(t\right).\label{eq:displacement-expansion}
\end{equation}

In a neighborhood of the manifold where $\bm{x}\approx\bm{y}$, these
have equations given by the chain rule for variable changes, noting
that we only consider smooth differentiable functions here:
\begin{equation}
\dot{\Delta}=\bm{J}\left(\bar{\bm{\phi}}+\Delta\right)\left[\bm{a}_{\lambda}\left(\bm{x}\left(\bar{\bm{\phi}}+\Delta\right)\right)+\bm{B}\left(\bm{x}\left(\bar{\bm{\phi}}+\Delta\right)\right)\bm{\xi}\right].
\end{equation}

\subsubsection{Fast variables}

Expanding the displacement $\bm{\Delta}$ according to Eq. (\ref{eq:displacement-expansion}),
we define
\begin{align}
a_{\perp}^{i} & =\bm{n}^{i}\cdot\bm{a}\nonumber \\
\bm{b}_{\perp}^{i} & =\bm{n}^{i}\cdot\bm{B}
\end{align}

Here $a_{\perp}^{i}$ is the $i$-th component of the normal drift,
and the inner-product notation refers to inner products in the original
Euclidean space. Each term $\bm{b}_{\perp}^{i}$ is a covariant vector
in the $s$-dimensional space of noise terms:
\begin{equation}
b_{\perp\sigma}^{i}=\sum_{j=1}^{n}n_{j}^{i}B_{\sigma}^{j}.
\end{equation}
This gives a locally valid stochastic equation for the fast variables,
where the inner product $\bm{b}_{\perp}^{i}\cdot\bm{\xi}=\sum_{\sigma}b_{\perp\sigma}^{i}\xi^{\sigma}$
is in the noise-vector space: 
\begin{equation}
\dot{\Delta}_{\perp}^{i}=a_{\perp}^{i}-\lambda\sum_{k=1}^{p}R^{ik}\Delta_{\perp}^{k}+\bm{b}_{\perp}^{i}\cdot\bm{\xi}.
\end{equation}

Here we define the constraint linear response matrix $\bm{R}$ as
\begin{equation}
R^{ik}=\sum_{j=1}^{p}\left(\boldsymbol{v}_{\perp}^{j}\cdot\bm{n}^{i}\right)\left(\boldsymbol{v}_{\perp}^{j}\cdot\bm{n}^{k}\right).
\end{equation}

One can rewrite $\bm{R}$ in a factored form as $\bm{R}=\bm{M}^{T}\bm{M}$,
where $M^{ij}=\boldsymbol{v}_{\perp}^{i}\cdot\bm{n}^{j}$ is the normal
projection matrix of (\ref{eq:normsl-projection-matrix}). From Eq
(\ref{eq:orthogonal-set-condition}), the diagonal elements of $\bm{R}$
are non-vanishing. As a result, $\bm{R}$ is a symmetric, positive
definite matrix, provided the constraints are non-singular. Since
it must have positive eigenvalues, the adiabatic limit is found by
taking the limit of $\lambda\rightarrow\infty$. This is equivalent
to setting $\dot{\Delta}_{j}^{\perp}=0$, due to the rapid equilibration
that occurs. 

We also assume that $\lambda T_{c}\gg0$ , with the result that the
``fast'' set $\Delta^{\perp}$ experience a large rate of change
defined by $\lambda$, and relaxes to equilibrium on a short time
scale of $1/\lambda$. Therefore, in the adiabatic limit, $\bm{\Delta}^{\perp}=0$,
and the system is constrained to the tangent space. 

\subsubsection{Slow variables}

By contrast, the slow set $\Delta^{\parallel}$ does not relax at
the same rate. The constraint terms vanish, since the Jacobian for
these variables is orthogonal to the constraint terms. Substituting
the resulting values into the ``slow'' equations therefore leads
to a simpler projected equation with fewer independent variables.
We define, for $\mu=1,\ldots m$ :
\begin{equation}
\dot{\Delta}^{\mu}=\sum_{j}J_{j}^{\mu}\left(a^{j}+B_{\sigma}^{j}\xi^{\sigma}\right)
\end{equation}

This equation can be written as an equation purely in terms of the
intrinsic manifold coordinates, as:
\begin{equation}
\dot{\Delta}^{\mu}=\alpha^{\mu}+\beta_{\sigma}^{\mu}\xi^{\sigma}
\end{equation}
with the definitions that:
\begin{align}
\alpha^{\mu} & =\sum_{j}J_{j}^{\mu}a^{j}\nonumber \\
\beta_{\sigma}^{\mu} & =\sum_{j}J_{j}^{\mu}B_{\sigma}^{j}.\label{eq:Intrinsic form}
\end{align}
This leads to an important question. Is the projected SDE an Ito or
Stratonovich equation? 

We assume that the original equation prior to projection has a finite
bandwidth. This means that ordinary calculus rules are applicable
at each stage. The broadband limit is then taken after adiabatic elimination,
subject to the restriction that $\lambda T_{c}\gg0$, which leads
to a Stratonovich interpretation. 

Defining $\bm{a}_{\parallel}\left(\bm{x}\right)=\mathcal{P}_{\bm{x}}^{\parallel}\left(\bm{a}\left(\bm{x}\right)\right)$
and $\bm{B}_{\parallel}\left(\bm{x}\right)=\mathcal{P}_{\bm{x}}^{\parallel}\left(\bm{B}\left(\bm{x}\right)\right)$,
the resulting Stratonovich equation can also be rewritten in a shorthand
form as:
\begin{equation}
\dot{\bm{x}}=\bm{a}_{\parallel}\left(\bm{x}\right)+\bm{B}_{\parallel}\left(\bm{x}\right)\bm{\xi}\left(t\right).\label{eq:Original-Strat-projected}
\end{equation}
Gardiner \citep{gardiner1984adiabatic} demonstrates that in two-dimensional
cases, then even if the original equation is in the broad-band limit,
adiabatic elimination will lead to the Stratonovich interpretation
of the resulting SDE, if there are no fast variables in the noise
coefficients. More generally, one should make a case-by-case analysis
to determine if there are additional stochastic correction terms.
The constraint equations are now automatically satisfied, provided
they are satisfied initially, since 
\begin{equation}
\dot{f}^{i}=\dot{\bm{x}}\cdot\nabla f^{i}=0.
\end{equation}

While these equations are correct in the continuous limit, it is important
to take measures to ensure the solutions remain on the manifold for
a finite time-step. To explain this, if we use the midpoint definition
of Eq. (\ref{eq:midpoint}), the discrete midpoint algorithm in Cartesian
coordinates can be written as:
\begin{equation}
\Delta\bm{x}=\bm{a}^{\parallel}\left(\bar{\bm{x}},\bar{t}\right)\Delta t+\bm{B}^{\parallel}\left(\bar{\bm{x}},\bar{t}\right)\Delta\bm{w}.
\end{equation}

Due to discretization error, the path may not remain on the manifold
with this method, as we show in numerical examples below. Since our
derivation requires that the constraints are satisfied exactly, we
will show that it is better to make an additional normal projection
to satisfy the constraints, so that the final algorithm reads: 

\begin{equation}
\Delta\bm{x}=\mathcal{P}^{\perp}\left\{ \bm{x}+\bm{a}^{\parallel}\left(\bar{\bm{x}},\bar{t}\right)\Delta t+\bm{B}^{\parallel}\left(\bar{\bm{x}},\bar{t}\right)\Delta\bm{w}\right\} -\bm{x}.
\end{equation}

\subsection{Intrinsic stochastic equations}

To summarize the preceding results, the projected coordinates $\bm{y}=\mathcal{P}\left(\bm{x}\right)$
are coordinates on an embedded manifold in Euclidean space $\mathbb{R}^{n}.$
A further transformation can be made to obtain intrinsic coordinates
on the manifold with lower dimensionality. Hence, to obtain $\bm{\phi}$,
we define

\begin{align}
\bm{\phi} & =\bm{\Phi}\left(\bm{y}\right).
\end{align}
One can make subsequent transformations on $\phi$ to obtain other
systems of coordinates on the manifold, such as defining spherical
polar coordinates with rotated polar directions, but these are all
isomorphic to the set of projected coordinates $\left\{ \mathcal{P}\left(\bm{x}\right)\right\} $. 

The diffusion and drift in (\ref{eq:Original-Strat-projected}) is
in a Stratonovich form \citep{Gardiner1997}, appropriate for many
physical problems. This is the wide-band limit of a physical noise,
and is valid in most cases of adiabatic elimination. As a result,
one can use ordinary calculus rules to transform this into an intrinsic
$n$-dimensional stochastic differential equation (SDE):
\begin{align}
\dot{\bm{\phi}} & =\bm{\alpha}\left(\bm{\phi}\right)+\bm{\beta}\left(\bm{\phi}\right)\bm{\xi},
\end{align}
where $\bm{\alpha}\left(\bm{\phi}\right)$, $\bm{\beta}\left(\bm{\phi}\right)$
are given by (\ref{eq:Intrinsic form}). The $s$-dimensional noise
vector $\bm{\xi}$ is Gaussian and delta-correlated, and follows Eq.
(\ref{eq:noise-correl}).

The corresponding diffusion or Fokker-Planck equation (FPE) \citep{Graham1977Covariant,Risken1996}
is for a probability density $P(\bm{\phi},t)$ in an $n$-dimensional
real vector space or manifold $\mathcal{M}$, where the intrinsic
coordinates are $\bm{\phi}$.  Our notation treats intrinsic phase-space
coordinates as contravariant vectors $\phi^{\mu}$, and derivatives
as covariant quantities $\partial_{\mu}\equiv\partial/\partial\phi^{\mu}$.
This leads to an FPE of form:
\begin{equation}
\frac{\partial P}{\partial t}=\left[-\partial_{\mu}\alpha^{\mu}(\bm{\phi})+\frac{1}{2}\partial_{\mu}\beta_{\sigma}^{\mu}(\bm{\phi})\partial_{\nu}\beta_{\sigma}^{\nu}(\bm{\phi})\,\right]P,\label{eq:StratonovichFPE}
\end{equation}
where we use the Einstein summation convention for repeated indices
$\mu=1,\ldots n$, $\sigma=1,\ldots s$, and the corresponding diffusion
matrix can also be regarded as a contravariant metric tensor \citep{Graham1977Covariant}.
To avoid confusion with the induced metric defined below, we use the
notation of a diffusion matrix $D^{\mu\nu}$,
\begin{equation}
D^{\mu\nu}(\bm{\phi})=\sum_{\sigma}\beta_{\sigma}^{\mu}(\bm{\phi})\beta_{\sigma}^{\nu}(\bm{\phi})\,.
\end{equation}
Growth restrictions on coefficients \citep{Arnold1992-stochastic}
are needed to guarantee that solutions exist. 

\subsection{Ito and covariant diffusion}

For comparison with previous work, there are several types of drift
term for stochastic equations, corresponding to other forms of differential
terms in the FPE and different types of stochastic calculus. The Stratonovich
drift $\bm{\alpha}$ used here is related to the Ito drift $\bm{\alpha}_{I}$
\citep{stratonovich1966new,Gardiner1997}, found in Ito stochastic
equations, by the mapping: 
\begin{equation}
\alpha_{I}^{\mu}=\alpha^{\mu}+\frac{1}{2}\sum_{\sigma}\beta_{\sigma}^{\nu}\partial_{\nu}\beta_{\sigma}^{\mu}.
\end{equation}
In the numerical examples, we take the simplest case that the original
manifold has a locally Euclidean metric and Euclidean diffusion $g_{ij}^{E}=\delta_{ij}$
on a tangent plane, although this is not the most general case. The
mapping to intrinsic coordinates defines an induced metric $g$, which
is in general a curved metric \citep{nash1956imbedding}, such that:
\begin{align}
g_{\mu\nu} & =\frac{\partial y^{i}}{\partial\phi^{\mu}}\frac{\partial y^{j}}{\partial\phi^{\nu}}\delta_{ij}.\label{eq:pullback}
\end{align}

We therefore consider an $m$-dimensional manifold equipped with a
metric $g_{\mu\nu}$, in an intrinsic coordinate system with coordinates
$\bm{\phi}.$ As an example, the usual Fokker-Planck equation governing
non-driven Euclidean diffusion on this manifold is given by
\begin{equation}
\frac{\partial\tilde{P}}{\partial t}=\frac{1}{2}D\tilde{\Delta}\tilde{P},
\end{equation}
where $\tilde{\Delta}=\nabla_{\mu}\nabla^{\mu}$ is the Laplace-Beltrami
operator, $\tilde{P}$ is the covariant probability density, and $D_{\text{}}$
is the covariant diffusion coefficient. This FPE is in a covariant
form, meaning that probability is conserved with respect to the measure
$\sqrt{g}d^{m}\bm{\phi}$, where $g=\det\left(g_{\mu\upsilon}\right).$
If one scales the covariant probability density by $\sqrt{g},$ the
result is a probability density $P$ that is conserved with respect
to $d^{m}\bm{\phi}$.

Inserting also the formula \citep{Risken1996,gustafsson1997diffusion}
, 
\begin{equation}
\tilde{\Delta}\tilde{P}=\frac{1}{\sqrt{g}}\partial_{\mu}(\sqrt{g}\partial^{\mu}\tilde{P}),
\end{equation}
one finds that
\begin{equation}
\frac{\partial P}{\partial t}=\frac{D}{2}\partial_{\mu}\left(\sqrt{g}g^{\mu\nu}\partial_{\nu}\frac{P}{\sqrt{g}}\right),
\end{equation}
with $g^{\mu\nu}=\left[g_{\mu\nu}\right]^{-1},$ and $P=\sqrt{g}\tilde{P}$.
Equivalently, using Einstein summation convention for repeated indices,
one obtains the self-adjoint diffusion equation,
\begin{equation}
\frac{\partial P}{\partial t}=-D\partial_{\mu}\left(\frac{1}{4}g^{\mu\nu}\partial_{\nu}\ln g\right)P+\frac{D}{2}\partial_{\mu}g^{\mu\nu}\partial_{\nu}P.
\end{equation}
In this expression the self-adjoint drift is 
\begin{equation}
\alpha_{a}^{\mu}=\frac{D}{4}g^{\mu\nu}\partial_{\nu}\ln g.\label{eq:SelfAdjointDrift}
\end{equation}
 This FPE may be transformed to the Stratonovich form, from which
one easily reads off an equivalent system of SDEs by introducing $Dg^{\mu\nu}=\sum_{\sigma}\beta_{\sigma}^{\mu}\beta_{\sigma}^{\nu}$,
and defining
\begin{equation}
\alpha_{\text{}}^{\mu}=\alpha_{a}^{\mu}+\frac{1}{2}\sum_{\sigma}\beta_{\sigma}^{\mu}\partial_{\nu}\beta_{\sigma}^{\nu},\label{eq:StratDrift}
\end{equation}
The corresponding projected Stratonovich SDE is then
\begin{equation}
\dot{\phi}^{\mu}=\alpha_{\text{}}^{\mu}+\beta_{\sigma}^{\mu}\xi^{\sigma}.
\end{equation}

However, the fact that the final diffusion matrix is proportional
to the induced metric is due to the form of the initial diffusion
matrix, which was chosen equal to the Euclidean metric tensor for
this example. This is not always the case, and the final diffusion
matrix has to be worked out accordingly.

\subsection{Hypersphere}

As an illustration, consider a hyper-spherical manifold with $\left|\bm{x}\right|^{2}=1$
and isotropic diffusion with $\bm{B}=\bm{I}$ in the original Euclidean
space. Writing Eq. (\ref{eq:Original-Strat}) as a difference equation
as in Eq. (\ref{eq:Finite-step-SDE}) leads to
\begin{equation}
\Delta\bm{x}=\Delta\bm{w}.
\end{equation}
where the Gaussian noise $\Delta\bm{w}$ integrated over an interval
$\Delta t$ is correlated according to Eq. (\ref{eq:Finite-step-noise-correlation}.

 Simply projecting this according to the constraint naively leads
to a projected equation: 

\begin{equation}
\Delta\bm{x}=\mathcal{P}^{\perp}\left(\bm{x}+\Delta\bm{w}\right)-\bm{x},
\end{equation}
where $\mathcal{P}^{\perp}\left(\bm{x}\right)$ projects the new
vector $\bm{y}=\bm{x}+\bm{\delta}$ normally onto the sphere, or more
generally onto a manifold $\mathcal{M}$. Instead of this approach,
using the adiabatic elimination combined with a potential, one might
choose that:
\begin{equation}
u\left(\bm{x}\right)=\frac{\lambda}{2}\left(\left|\bm{x}\right|^{2}-1\right)^{2}.\label{eq:spherical-constraint-potential-1}
\end{equation}
On the projected manifold, the restoring drift towards the surface
is normal, since: 
\begin{align}
\nabla u\left(\bm{y}\right) & =2\lambda\bm{y}\left(\left|\bm{y}\right|^{2}-1\right).
\end{align}
For a tangential projection $\mathcal{P}_{\bm{x}}^{\parallel}$ on
the manifold, the resulting stochastic equation is then:
\begin{equation}
\dot{\bm{x}}=\mathcal{P}_{\bm{x}}^{\parallel}\bm{\xi}=\bm{\xi}_{\bm{x}}^{\parallel}.\label{eq:noise-projection}
\end{equation}

From the results above, this is a Stratonovich equation \citep{stratonovich1966new}.
The corresponding intrinsic equations for the $3$-dimensional case
(see Appendix) have an unphysical singularity at the poles, which
is eliminated using projections:
\begin{align}
\dot{\theta} & =\frac{1}{2}\cot\theta+\xi^{\theta}\nonumber \\
\dot{\phi} & =\frac{\xi^{\phi}}{\sin\theta}.
\end{align}

\section{Stochastic projection algorithms \label{sec:Stochastic-projection-algorithms}}

We first summarize the numerical algorithms for tangential and normal
projections, and then explain how these are used in projected stochastic
equations. In some cases \citep{armstrong2019optimal}, a tangential
projection of the coefficients of a Stratonovich SDE is used. However,
this doesn't guarantee the solution is on the manifold. The difficulty
with tangential methods is that global errors can accumulate at finite
step-size. Hence, the final result can move arbitrarily far way from
the desired manifold. Other proposed methods that use restrictions
to a manifold are often limited to two-dimensional surfaces in three-dimensional
spaces \citep{yang2019simulation}. 

Another approach commonly applied to physics problems \citep{holyst1999diffusion},
is to employ a tangentially projected Ito-Euler step followed by normal
projections to remain on the manifold. Similar techniques are used
in molecular dynamics applications to compute free energies \citep{ciccotti2008projection,Lelievre2010free}. 

Compared to a purely tangential approach, this has the advantage that
the final normal projection keeps the solution on the manifold. However,
Ito stochastic equations are not always equivalent to Stratonovich
stochastic equations, unless appropriate corrections are employed. 

The hybrid approach described here is to use a stochastic method of
known convergence properties to solve the tangentially projected Stratonovich
equations, (\ref{eq:Original-Strat-projected}), followed by normal
projection to the manifold after each step. This combines the best
features of both the tangential and normal projection approaches.
There are also stochastic methods of higher orders available \citep{Kloeden1992}.
We refer the reader to recent work in this area \citep{laurent2021order}.

 Our approach is illustrated with a hybrid midpoint projection algorithm
that combines a midpoint method \citep{Drummond1990} instead of an
Euler step, together with normal projection. The manifold is given
by the constraint equations $f^{i}\left(\bar{\bm{x}}\right)=0$, where
$\bar{\bm{x}}$ is defined as a midpoint for each tangential step.
We assume the existence of an orthogonalization algorithm that generates
an orthonormal set of vectors , $S^{\perp}=\left\{ \bm{n}^{1},\ldots\bm{n}^{p}\right\} $,
from a set of gradient vectors $S=\left\{ \boldsymbol{v}_{\perp}^{j},\ldots\boldsymbol{v}_{\perp}^{p}\right\} $,
using one of the known standard techniques \citep{trefethen1997numerical,lowdin1950non}. 

The three algorithms compared here all use the notation that $\bm{x}_{0},\bm{x}_{1}$
are the initial and final locations of a step in time, $\Delta t$
is the step-size, and $\Delta\bm{w}$ are random Gaussian variables.
They have correlations that correspond to a discretized delta-function
in time, as described in Eq (\ref{eq:Finite-step-noise-correlation}).

In the next section, we give examples of the use of these algorithms
by comparing them either with exact results or with well-converged,
high-accuracy simulations, for a wide range of dimensions and different
constraint types, namely:
\begin{itemize}
\item Kubo oscillators
\item catenoids
\item hyperspheres 
\item spheroids
\item hyperboloids
\item polynomial surfaces
\end{itemize}

\subsection{Error properties}

There are three main errors when solving projected stochastic equations
numerically. The first is caused by the finite sample number, $N_{s}$.
This means that probabilities or averages over the ensemble of paths
have a sampling error of order $1/\sqrt{N_{s}}$, which depends on
the noise term. Such errors can be reduced further \citep{Kloeden1992,opanchuk2016parallel,kiesewetter2017algorithms},
but this involves additional complexity. Reducing sampling errors
requires fast algorithms to increase $N_{s}.$

The second is the step-size or discretization error, caused by the
finite size of the time-step $\Delta t$. Because noise terms have
fluctuations that scale as $\Delta w\sim\sqrt{\Delta t}$, these errors
requires different numerical algorithms to achieve a given convergence
order than with ordinary differential equations. 

Due to the differences between drift and noise, one can usefully specify
the single-step error order as $\left[\Delta t^{n},\Delta W^{m}\right]$
\citep{Drummond1990}, where the first term quantifies the zero-noise
limit, and the second term the noise error, scaling as $\Delta t^{m/2}$.
The global error after finite time $T$ usually scales as $\left[\Delta t^{n-1},\Delta W^{m-1}\right]$
\citep{mil1975approximate}, in the small step-size limit.

There is a third type of error, which is the global error in the projection
constraint. We regard this as a distinct error, depending on the detailed
constraint equations.

One can also distinguish between weak and strong convergence errors,
where weak convergence measures errors in the ensemble averaged probabilities
$P\left(\bm{x},t\right)$ and its moments, while strong convergence
measures errors of a given random trajectory $\bm{x}\left(t\right)$
conditioned on a specific noise $\bm{\xi}\left(t\right)$.

All stochastic methods using independent trajectories have the same
scaling law for sampling errors. The local discretization and global
constraint errors depend on the method used, and are treated in greater
detail below. 

We focus on relatively low order, high performance methods, as they
are often most useful in applications. The reason for this is that
the total error combines discretization and sampling errors. Reducing
the discretization error using high-order methods is not as useful
as with ordinary differential equations, since the increased algorithmic
complexity and need for small step-sizes for convergence may lead
to fewer samples and increased sampling error.

The algorithms assume a common toolbox of tangential and normal projection
methods, denoted as $tang\left(\bm{\delta}\left|\bm{x}_{0}\right.\right)$
and $norm\left(\bm{x}\right)$ respectively, and defined below in
detail.

\subsection{Combined Euler projection (cEP) \label{subsec:xEPtn}}

Holyst \citep{holyst1999diffusion}, describes a combined Euler projected
method with pure diffusion. The original proposal had no drift vector,
$\bm{a}=0$, and was restricted to dimensions $n=3$ and $m=2.$ Their
approach uses an Ito-Euler step in time, with a combined tangential
projection of the noise at the initial point, and a normal projection
at the end of the step in time. Here we describe this approach more
generally by including a drift term and a constant noise matrix $\bm{B}$.
We note that this Ito-type algorithm requires additional corrections
if the diffusion is space-dependent. Similar techniques are used in
molecular dynamics \citep{ciccotti2008projection,Lelievre2010free}. 

In the notation of Eq ( \ref{eq:Finite-step-SDE}), the algorithm
is defined for a step starting at $\bm{x}_{0}$, assuming that $\bm{a}$
is the drift and $\bm{B}$ is the noise matrix, as:

\begin{equation}
\Delta\bm{x}=\mathcal{P}^{\perp}\left(\bm{x}+\left(\bm{a}^{\parallel}\left(\bm{x}_{0}\right)\Delta t+\bm{B}^{\parallel}\left(\bm{x}_{0}\right)\cdot\Delta\bm{w}\right)\right)-\bm{x}.\label{eq:Finite-step-SDE-1}
\end{equation}

The projected algorithm is described in detail in the box below: 
\begin{center}
\doublebox{\begin{minipage}[t]{0.9\columnwidth}%
\begin{enumerate}
\item Evaluate a step, $\bm{\delta}=\bm{a}\left(\bm{x}_{0}\right)\Delta t+\bm{B}\left(\bm{x}_{0}\right)\cdot\Delta\bm{w}.$
\item Tangentially project: $\bm{\delta}_{\parallel}=tang\left(\bm{\delta}\left|\bm{x}_{0}\right.\right)$
\item Estimate an intermediate point: $\bm{x}'=\bm{x}_{0}+\bm{\delta}_{\parallel}$
\item Return $\bm{x}_{1}=norm\left(\bm{x}'\right)$
\end{enumerate}
\end{minipage}}
\par\end{center}

\subsection{Tangential midpoint projection (tMP) \label{subsec:xGPt}}

We next consider purely tangential projections using a midpoint projection
algorithm which can be expressed symbolically as:

\begin{equation}
\Delta\bm{x}=\bm{a}^{\parallel}\left(\bar{\bm{x}},\bar{t}\right)\Delta t+\bm{B}^{\parallel}\left(\bar{\bm{x}},\bar{t}\right)\Delta\bm{w}.\label{eq:Finite-step-SDE-1-1}
\end{equation}

This method is identical to a known implicit midpoint method \citep{Drummond1990},
except with a tangential drift and noise. It is expected to handle
stiffness in the equations better than explicit Euler-based methods
\citep{Kloeden1992}, and does not require knowledge of the derivatives. 

From previous analysis, it has a local accuracy that depends on the
properties of the noise coefficient $\bm{B}^{\parallel}$. For general
noise coefficients, the strong local error is of order $\left[\Delta t^{3},\Delta W^{2}\right]$,
although this can be improved if the noise is commutative. 

The relevant quantity for calculating probabilities and statistics
at a finite time is the weak global error at finite time, which is
of order $\left[\Delta t^{2},\Delta W^{2}\right]$ \citep{Drummond1990,Kloeden1992}
. 

Although it converges with the given accuracy in the limit of small
step-size, these conditions do not guarantee that a path stays on
the manifold. This can result in substantial global errors at finite
step-size.We note that a Stratonovich-type midpoint algorithm evaluates
derivatives at the midpoint \citep{Drummond1990}: 
\begin{center}
\doublebox{\begin{minipage}[t]{0.9\columnwidth}%
\begin{enumerate}
\item Set $m=0$ and $\bar{\bm{x}}^{\left(0\right)}=\bm{x}_{0}.$
\item $\bm{\delta}=\left[\bm{a}\left(\bar{\bm{x}}^{\left(m\right)}\right)\Delta t+\bm{B}\left(\bar{\bm{x}}^{\left(m\right)}\right)\cdot\Delta\bm{w}\right]/2$
\item Tangentially project: $\bm{\delta}_{\parallel}=tang\left(\bm{\delta}\left|\bar{\bm{x}}^{\left(m\right)}\right.\right)$
\item Estimate midpoint: $\bar{\bm{x}}^{\left(m+1\right)}=\bm{x}_{0}+\bm{\delta}_{\parallel}$
\item If $m<iters$, $m\rightarrow m+1$, go to $\left(2\right)$
\item Return $\bm{x}_{1}=\bm{x}_{0}+2\bm{\delta}_{\parallel}$
\end{enumerate}
\end{minipage}}
\par\end{center}

This method has lower local errors than the Euler method, but as is
not constrained to say on the manifold, there is an increasing constraint
error with time. This problem is removed by using an additional normal
projection. 

\subsection{Combined midpoint projection (cMP)\label{subsec:xGP}}

Combined tangential and normal projection leads to an algorithm for
the projected stochastic differential equation which has both good
accuracy and long-term constraint stability. 

It can be written symbolically as:
\begin{equation}
\Delta\bm{x}=\mathcal{P}^{\perp}\left(\bm{x}+\bm{a}^{\parallel}\left(\bar{\bm{x}},\bar{t}\right)\Delta t+\bm{B}^{\parallel}\left(\bar{\bm{x}},\bar{t}\right)\Delta\bm{w}.\right)-\bm{x}
\end{equation}
 Both the midpoint solution and the final normal projection are obtained
with fixed-point iteration.
\begin{center}
\doublebox{\begin{minipage}[t]{0.9\columnwidth}%
\begin{enumerate}
\item Set $m=0$ and $\bar{\bm{x}}^{\left(0\right)}=\bm{x}_{0}.$
\item $\bm{\delta}=\left[\bm{a}\left(\bar{\bm{x}}^{\left(m\right)}\right)\Delta t+\bm{B}\left(\bar{\bm{x}}^{\left(m\right)}\right)\cdot\Delta\bm{w}\right]/2$
\item Tangentially project: $\bm{\delta}_{\parallel}=tang\left(\bm{\delta}\left|\bar{\bm{x}}^{\left(m\right)}\right.\right)$
\item Estimate midpoint: $\bar{\bm{x}}^{\left(m+1\right)}=\bm{x}_{0}+\bm{\delta}_{\parallel}$
\item If $m<iters$, $m\rightarrow m+1$, go to $\left(2\right)$
\item Return $\bm{x}_{1}=norm\left(\bm{x}_{0}+2\bm{\delta}_{\parallel}\right)$
\end{enumerate}
\end{minipage}}
\par\end{center}

This method uses a tangential projection for each midpoint iteration,
together with a final normal projection to ensure that the resulting
step remains on the manifold. The reason for this choice is that the
stochastic departure from the manifold is significant if allowed to
propagate. The final normal  projection removes this additional error.

The final projection step, if convergent, can only reduce the global
error, as it constrains trajectories to the manifold with even greater
accuracy. We have verified this numerically with exact examples, given
below. Since the constraint is satisfied through a final normal projection,
one expects that the global constraint error will be very small, and
in some examples it reaches the limit of roundoff errors, about $10^{-16}$
with IEEE digital arithmetic.

It is important to ensure that the projected path remains on the manifold,
since error-propagation often may cause a long term drift off the
manifold and rapid error growth. There are two parts to the calculation,
which in principle would require different numbers of iterations.
In the examples we use three iterations for each part. These fixed
point iterations may not converge at large step-sizes, so it is essential
to compare results at different step-sizes to check this.

\subsection{Quantitative error comparisons}

Although the three methods compared here are first order, they have
quite different errors in practical terms. To understand this, consider
a unit circular projection in the simplest case of the driven Kubo
oscillator, treated below, equivalent to $\dot{z}=iz\left(\omega\left(t\right)+b\xi\left(t\right)\right)$,
where $z=x^{1}+ix^{2}$, and $\left|z\right|=1$. Defining $z=\exp\left(i\theta\right)$
gives a readily soluble exact intrinsic equation with:
\begin{equation}
\dot{\theta}=\omega\left(t\right)+b\xi\left(t\right).
\end{equation}

This has the solution that: 
\begin{equation}
\theta\left(t\right)=\theta\left(0\right)+\int_{0}^{t}\omega\left(\tau\right)d\tau+b\Delta w\left(t\right),
\end{equation}

 where $\left\langle \Delta w^{2}\left(t\right)\right\rangle =t$.
Comparing exact results $\left\langle O\left(t\right)\right\rangle $
for an observable $O$, with numerical results $\left\langle O_{A}\left(t\right)\right\rangle $
using an algorithm $A$ gives a truncation error defined as: 
\begin{equation}
e\left(\left\langle O_{A}\left(t\right)\right\rangle \right)=\left\langle O_{A}\left(t\right)\right\rangle -\left\langle O\left(t\right)\right\rangle .
\end{equation}

For the normally projected methods, angular errors are the only significant
error, but for the methods without normal projection there is an additional
radial or projection error, which grows in time as well. 

\subsubsection{Euler projections}

For Euler projected methods starting at $t_{0}$, with $\omega_{0}=\omega\left(t_{0}\right)$,
the angular change $\Delta\theta$ for a finite step $\Delta=\omega_{0}\Delta t+b\Delta w$
has a projection of:
\begin{equation}
\Delta\theta_{E}=\arctan\left(\Delta\right)\approx\Delta-\frac{1}{3}\Delta{}^{3}+O\left(\Delta{}^{5}\right).
\end{equation}
Since $\left\langle \Delta w^{4}\right\rangle =3\Delta t^{2}$, this
leads to mean-squared displacements for pure noise evolution, of 
\begin{equation}
\left\langle \Delta\theta_{E}^{2}\right\rangle =b^{2}\Delta t-2b^{4}\Delta t^{2}+O\left(\Delta t^{3}\right).
\end{equation}

Starting from $\theta_{0}=0$, the global $\left\langle \theta_{E}^{2}\right\rangle $
error in the pure diffusion limit is therefore:
\begin{equation}
\left|e\left(\left\langle \theta_{E}^{2}\right\rangle \right)\right|\le2tb^{4}\Delta t.\label{eq:global-error-Euler}
\end{equation}
This leads to first order convergence, with global errors scaling
as $\Delta t$, just as in the pure drift case.

\subsubsection{Midpoint projections}

By comparison, the midpoint projection leads to a much smaller local
error of $\Delta{}^{3}/24$
\begin{equation}
\Delta\theta_{MP}=2\arcsin\left(\Delta/2\right)\approx\Delta+\frac{1}{24}\Delta{}^{3}.
\end{equation}
 For pure noise, this leads to mean-squared displacements with local
errors smaller by a factor of $8$. However, this error estimate is
only valid with a step that commences on the manifold, in this case
with $\left|z\right|=1$. The projection error can grow in time unless
a second normal projection is used. Starting from $\theta_{0}=0$,
the global $\left\langle \theta_{MP}^{2}\right\rangle $ error with
combined projections in the pure diffusion limit is therefore:
\begin{equation}
\left|e\left(\left\langle \theta_{MP}^{2}\right\rangle \right)\right|\le\frac{1}{4}tb^{4}\Delta t.\label{eq:Global-error-MP}
\end{equation}

This explains the much lower global errors observed numerically for
cMP compared to cEP in the examples below. Including drift terms gives
even greater improvements, since the method has a second order convergence
in the noise-free limit with $b=0$. 

\subsubsection{Higher-order methods}

Since the original equation was in a Stratonovich form, any higher-order
method can also be used \citep{Kloeden1992}, provided it is designed
for a Stratonovich process and includes a tangential projection in
each evaluation of the derivative. Based on the numerical experiments
carried out here, one should use a normal projection to prevent drift
off the manifold. Related methods have been used in free-energy sampling
\citep{laurent2021order}.

\subsection{Projection algorithms}

We use numerical projection algorithms defined as follows, where $\bm{\delta}$
is an estimated derivative step and $\bar{\bm{x}}$ is an intermediate
estimate:

\subsubsection{Tangential projection of $\bm{\delta}$ at $\bar{\bm{x}}$: $\bm{\delta}_{\parallel}=tang\left(\bm{\delta}\left|\bar{\bm{x}}\right.\right)$}
\begin{center}
\doublebox{\begin{minipage}[t]{0.9\columnwidth}%
\begin{enumerate}
\item Obtain normal vectors $\boldsymbol{v}_{\perp}^{j}$ at $\bar{\bm{x}}$
\item Calculate orthonormal gradient vectors $\bm{n}^{j}$ from $\boldsymbol{v}_{\perp}^{j}$
\item Return $\bm{\delta}_{\parallel}=\bm{\delta}-\sum_{j=1}^{m}\bm{n}^{j}\left(\bm{\delta}\cdot\bm{n}^{j}\right)$
\end{enumerate}
\end{minipage}}
\par\end{center}

\subsubsection{Normal projection at $\bar{\bm{x}}$: $\bm{x}=norm\left(\bar{\bm{x}}\right)$}
\begin{center}
\doublebox{\begin{minipage}[t]{0.9\columnwidth}%
\begin{enumerate}
\item Set $m=0$ and $\bm{x}^{(0)}=\bar{\bm{x}}.$
\item Obtain normal vectors $\boldsymbol{v}_{\perp}^{j}$ and $\bm{n}^{j}$
at $\bm{x}^{(m)}$
\item Evaluate $M^{ij}=\bm{v}_{\perp}^{i}\cdot\boldsymbol{n}^{j}$
\item $\bm{x}^{(m+1)}=\bm{x}_{\parallel}^{(m)}-\sum_{i,j}\bm{n}^{i}\left[\bm{M}\right]_{ij}^{-1}f^{j}\left(\bm{x}^{(m)}\right)$
\item If $m<max$, $m\rightarrow m+1$, go to $\left(2\right)$
\item Return $\bm{x}=\bm{x}^{(m+1)}$
\end{enumerate}
\end{minipage}}
\par\end{center}

\section{Numerical examples \label{sec:Numerical-examples}}

In this section we give numerical examples for the algorithms. Testing
convergence by reducing the step-size does not always verify that
the converged result is free of errors. Therefore we use exact results
and comparisons with intrinsic methods to test and verify the projection
algorithms. 

An alternative method for comparisons is to check against steady-state
distributions. This is restricted to long-time comparisons in cases
where analytic solutions exist \citep{laurent2021order}. Exact comparisons
cannot be carried out in general, since intrinsic coordinates and
exact results are not always known.

Our numerical comparisons use a public domain stochastic toolbox \citep{kiesewetter2016xspde}
to dynamically compare different curved space diffusion SDE projection
methods with an intrinsic method. We consider circular, hyper-spherical,
spheroidal, catenoidal and hyperboloidal surfaces, as well as a higher
dimensional polynomial surface.

These examples have positive and negative curvatures which vary in
space, as well as higher dimensions and non-quadratic constraints,
to cover a variety of conditions. Errors were calculated by comparison
either with intrinsic simulations or with a case using a very small
step-size. 

Details of the Stratonovich equations used for the intrinsic calculations
are given in the Appendix. Errors in the reference calculations were
reduced as follows:
\begin{itemize}
\item The sampling errors were reduced by using $10^{7}$ trajectories,
which gives sampling errors in the range $\sim10^{-4}-10^{-3}$ in
the means, and in all cases smaller than step-size errors.
\item Step-size reductions of a factor of $5$ were used for comparisons
if there was no exact solution. These had five times smaller errors
than the tests, so these errors were negligible.
\item An implicit midpoint algorithm \citep{Drummond1990} was used for
intrinsic coordinates, otherwise the most accurate projected method
was used for comparisons.
\end{itemize}
Step-sizes of $\Delta t=0.1$ and smaller were used to investigate
scaling, with unit diffusion. Timings of the three methods were similar,
in the ratio $0.5:0.7:1$ for tMP:cEP:cMP algorithms. This is because
the second projection requires fixed-point iteration, so combined
methods are slower than one using tangential methods only, as one
would expect.

In all figures, the two solid lines, where visible, are the upper
and lower $\pm\sigma$ bounds from sampling errors in the projection
algorithms with $10^{7}$ trajectories, and the error bars indicate
the sampling errors in the reference simulations.

\subsection{Kubo oscillator}

The simplest example is the projection of Euclidean noise onto a circle,
giving the Kubo oscillator which is widely used as a model of microscopic
noise in solid-state physics \citep{Anderson1954mathematical,kubo1954note,chaudhuri2009microscopic}.
If $\bm{x}=\left(x,y\right)$ and $\mathcal{P}_{\bm{x}}^{\parallel}$
projects onto a unit circle with $f\left(\bm{x}\right)=x^{2}+y^{2}-1$,
this can be written as:
\begin{equation}
\dot{\bm{x}}=\mathcal{P}_{\bm{x}}^{\parallel}\left(\bm{a}+\bm{\xi}\right),
\end{equation}
where $\bm{a}$ is a drift term. Defining $z=x+iy$, the projected
form has a complex Stratonovich equation, $\dot{z}=iz\left(\omega\left(t\right)+b\xi\left(t\right)\right)$.
The corresponding intrinsic form is $\dot{\theta}=\omega\left(t\right)+b\xi\left(t\right)$,
where $z=\exp\left(i\theta\right)$.

This has exact solutions for all complex moments:
\begin{equation}
\left\langle \left[z\left(t\right)\right]{}^{m}\right\rangle =\left\langle \left[z\left(0\right)\right]{}^{m}\right\rangle e^{\left(im\int_{0}^{t}\omega\left(\tau\right)d\tau-m^{2}b^{2}t/2\right)}.
\end{equation}
For example, if $z\left(0\right)=b=1$ and $\omega\left(t\right)=\omega_{0}t$
to give a definite case, then:
\begin{align}
\left\langle x\left(t\right)\right\rangle  & =e^{-t/2}\cos\left(\omega_{0}t^{2}/2\right).
\end{align}

\begin{figure}[H]
\centering{}\includegraphics[width=0.75\columnwidth]{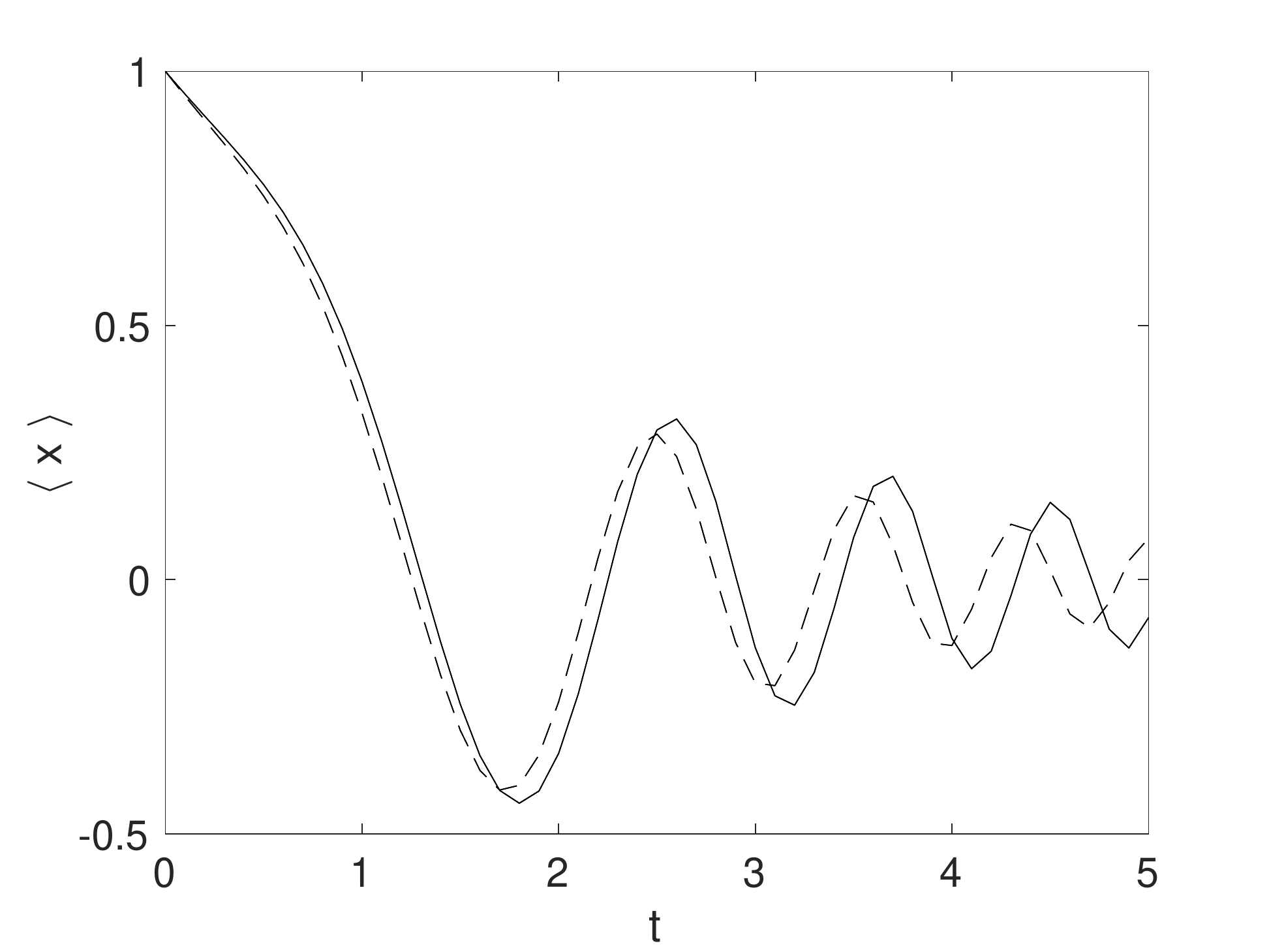}\\
\includegraphics[width=0.75\columnwidth]{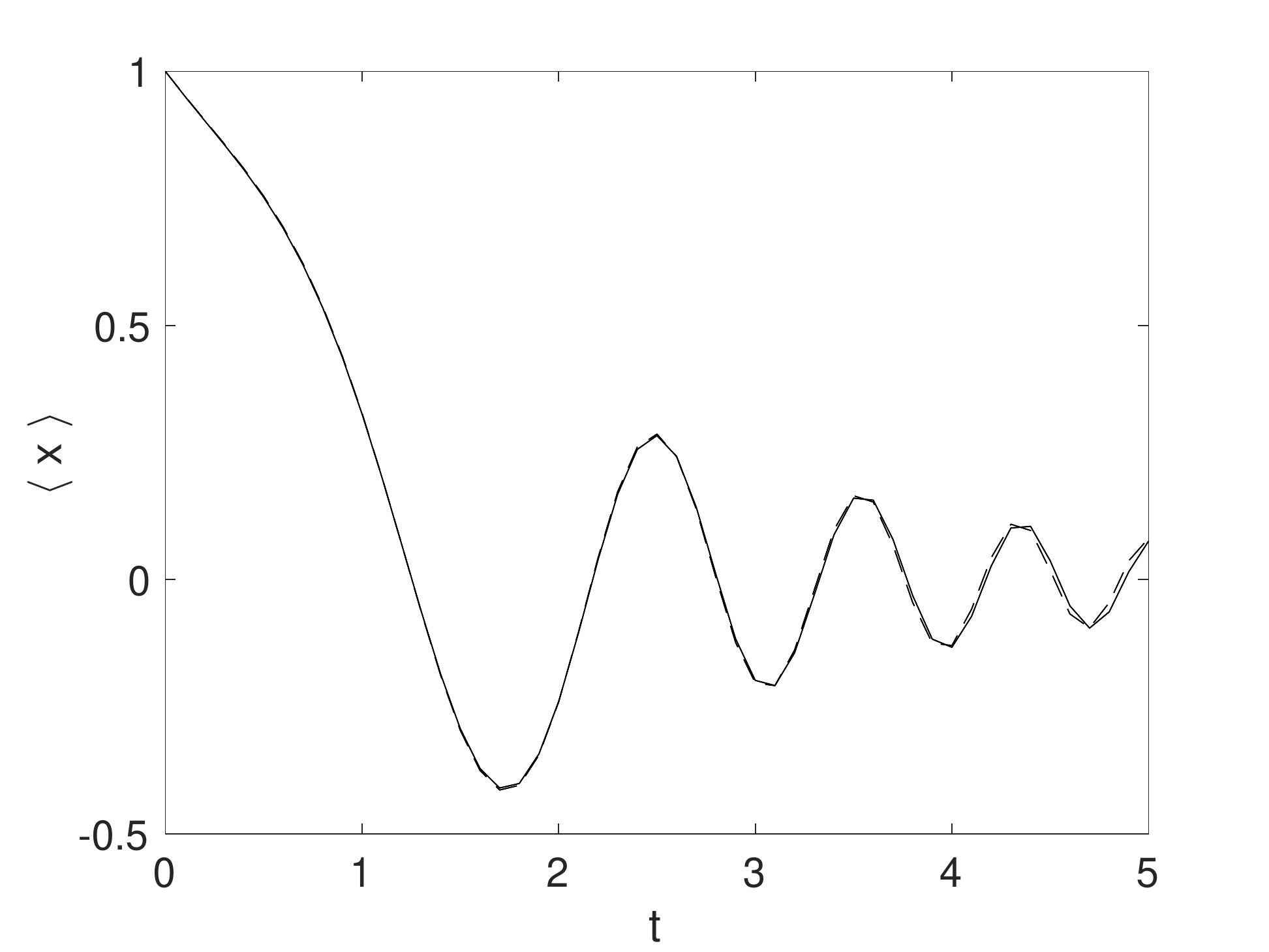}\caption{Comparing projection algorithms for the Kubo oscillator case with
a step-size of $0.05$. The solid lines are simulations, dashed lines
are the exact results. Top: combined Euler projection method (cEP).
Bottom, combined midpoint projection method (cMP). \label{fig:Kubo}}
\end{figure}

Table (\ref{tab:Comparison-of-error-Kubo}) gives comparative errors
using the different algorithms, for the case $b=1$ and $\omega_{0}=2.5$.
A comparative graph of mean values is given in Fig(\ref{fig:Kubo}).

For a step-size $\Delta t=0.05$, the combined midpoint method has
$\sim10$ times lower maximum errors in $\left\langle x\right\rangle $
than the Euler method, and over $10^{12}$ times lower maximum projection
errors than the tangential method. The quadratic error reduction in
the combined midpoint method occurs because the step-size error is
dominated by the drift term. In this limit, the midpoint method has
second-order global error convergence. Sampling errors are of order
$2\times10^{-4}$.

\begin{table}[H]
\centering{}%
\begin{tabular}{|c|c|c|c|c|}
\hline 
Function & Step & cEP & tMP  & \textbf{cMP }\tabularnewline
\hline 
\hline 
$\left\langle x\right\rangle $  & $0.1$ & $0.38$ & \textbf{$0.11$} & \textbf{$0.11$}\tabularnewline
\hline 
 & $0.05$ & $0.19$ & \textbf{$2\times10^{-2}$} & \textbf{$2\times10^{-2}$}\tabularnewline
\hline 
$\left\langle \left|f\right|\right\rangle $ & $0.1$ & $3.7\times10^{-5}$ & $0.3$ & \textbf{$1.2\times10^{-5}$}\tabularnewline
\hline 
 & $0.05$ & $1.8\times10^{-7}$ & $0.03$ & \textbf{$6\times10^{-16}$}\tabularnewline
\hline 
\end{tabular}\caption{Comparison of maximum error of combined Euler projection (cEP), tangential
midpoint projection (tMP), and combined midpoint projections (cMP)
for Kubo oscillator, using time-steps of $0.1$ and 0.05, with $10^{7}$
parallel trajectories, $t_{max}=5$ and $\omega_{0}=2.5$. \label{tab:Comparison-of-error-Kubo}}
\end{table}

\subsection{Catenoid}

The catenoid is an hyperbola rotated around an axis, and is a surface
of zero curvature. Despite the fact that the geometry is non-planar,
it has linear growth in diffusion distance with time, just as a planar
surface has. This is a minimal surface, defined by the constraint
that \citep{berger2012differential}:

\noindent 
\begin{equation}
x^{2}+y^{2}-\sinh^{2}z=1.
\end{equation}
A convenient set of intrinsic coordinates is \citep{berger2012differential}:
\begin{align}
x & =\cosh v\cos\theta,\nonumber \\
y & =\cosh v\sin\theta,\nonumber \\
z & =v,
\end{align}

For unit diffusion, so $D=B=I$, the exact solution for the average
diffusion distance is \citep{castro2014intrinsic}:

\begin{equation}
\left\langle \left|\bm{x}-\bm{x}_{0}\right|^{2}\right\rangle =2t.
\end{equation}

A maximum error comparison is shown in Table (\ref{tab:Comparison-of-error-catenoid}),
using time-steps of $0.1$ and 0.05, with a duration of $t_{max}=5$
and $10^{7}$ parallel trajectories. Sampling errors in distance are
of order $3\times10^{-3}$. 
\begin{table}[H]
\centering{}%
\begin{tabular}{|c|c|c|c|c|}
\hline 
Function & $\Delta t$ & cEP & tMP & \textbf{cMP}\tabularnewline
\hline 
\hline 
$\left\langle R^{2}\right\rangle $ & $0.1$ & $0.45$ & $7.8\times10^{-2}$ & $5.7\times10^{-2}$\tabularnewline
\hline 
 & $0.05$ & $0.24$ & $4.4\times10^{-2}$ & $3.3\times10^{-2}$\tabularnewline
\hline 
$\left\langle \left|f\right|\right\rangle $ & $0.1$ & $7.6\times10^{-8}$ & $7.3\times10^{-2}$ & $3.2\times10^{-13}$\tabularnewline
\hline 
 & $0.05$ & $2.2\times10^{-10}$ & $4\times10^{-2}$ & $1.6\times10^{-15}$\tabularnewline
\hline 
\end{tabular}\caption{Comparison of maximum error of combined Euler projection (cEP), tangential
midpoint projection (tMP), and combined midpoint projection (cMP)
methods for diffusion on a catenoid. \label{tab:Comparison-of-error-catenoid}}
\end{table}

For a step-size of $0.05,$ the  advantage of the cMP method is a
factor of more than $8$ reduction in mean square distance error compared
to combined Euler projection, and more than $10^{12}$ reduction in
constraint error compared to tangential midpoint projection.

\subsection{Spheroid}

An ellipsoid is defined by the constraint equation \citep{gray2017modern}
\begin{equation}
x^{2}/a^{2}+y^{2}/b^{2}+z^{2}/c^{2}-1=0.
\end{equation}
Setting $a=b=1$ gives a spheroid, which is a distortion of the unit
sphere along one axis. The metric tensor on the spheroid simplifies
greatly compared to the general ellipsoid, and we restrict ourselves
to this shape. In order to test the projection algorithms, we compare
the results of the projected SDE to a direct simulation of diffusion
on the spheroid using the intrinsic coordinates of:
\begin{align}
x & =\sin\theta\cos\phi,\nonumber \\
y & =\sin\theta\sin\phi,\nonumber \\
z & =c\cos\theta,
\end{align}
with $0\leq\theta\leq\pi$ and $0\leq\phi<2\pi.$ The equations used
for this are derived in the Appendix, and are:

\begin{align}
\dot{\theta} & =-\frac{\cot\theta}{(c^{2}-1)\cos2\theta-(1+c^{2})}+\frac{\xi^{\theta}}{\sqrt{\cos^{2}\theta+c^{2}\sin^{2}\theta}}\nonumber \\
\dot{\phi} & =\frac{\xi^{\phi}}{\sin\theta}.
\end{align}

\begin{figure}[H]
\begin{centering}
\includegraphics[width=0.75\columnwidth]{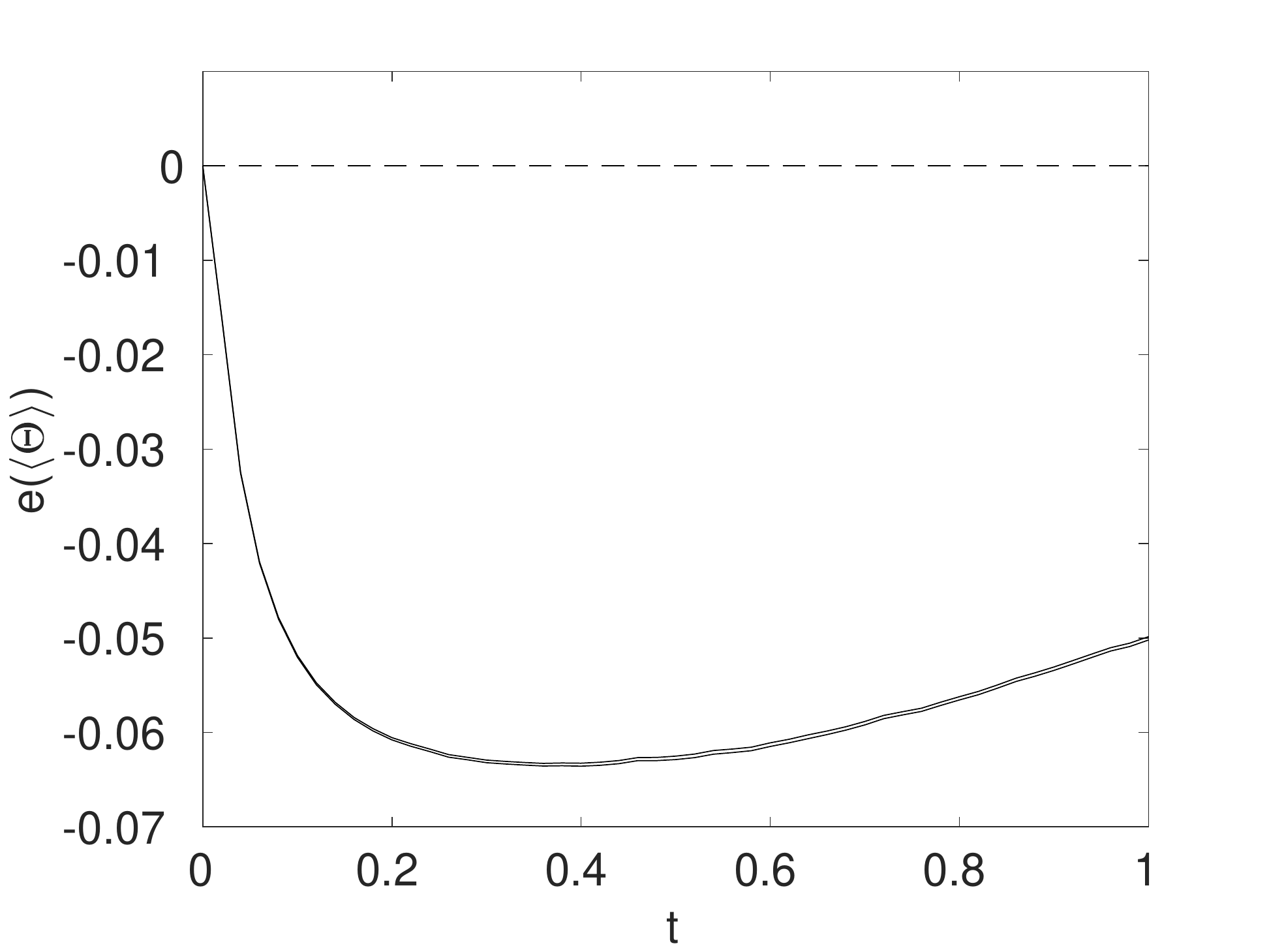}
\par\end{centering}
\centering{}\caption{cEP mean great circle distance error, $e\left(\left\langle \Theta\right\rangle \right)$
for diffusion on a spheroid. \label{fig:Ellipsoid-hEP}}
\end{figure}

\begin{figure}[H]
\begin{centering}
\includegraphics[width=0.75\columnwidth]{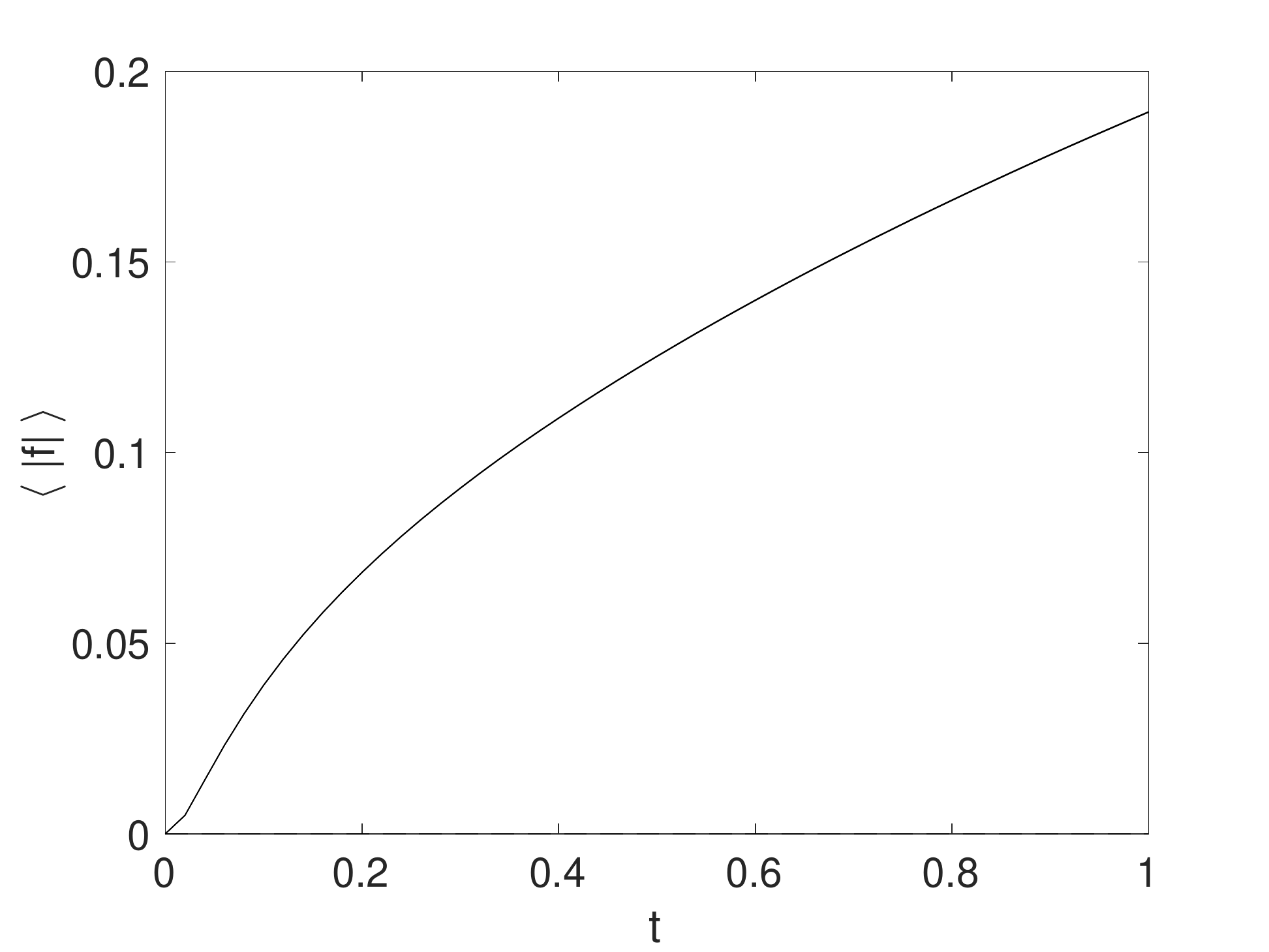}
\par\end{centering}
\caption{(Constraint error of the tMP algorithm for diffusion on a spheroid.
\label{fig:Ellipsoid-MP}}
\end{figure}

\begin{figure}[H]
\centering{}\includegraphics[width=0.75\columnwidth]{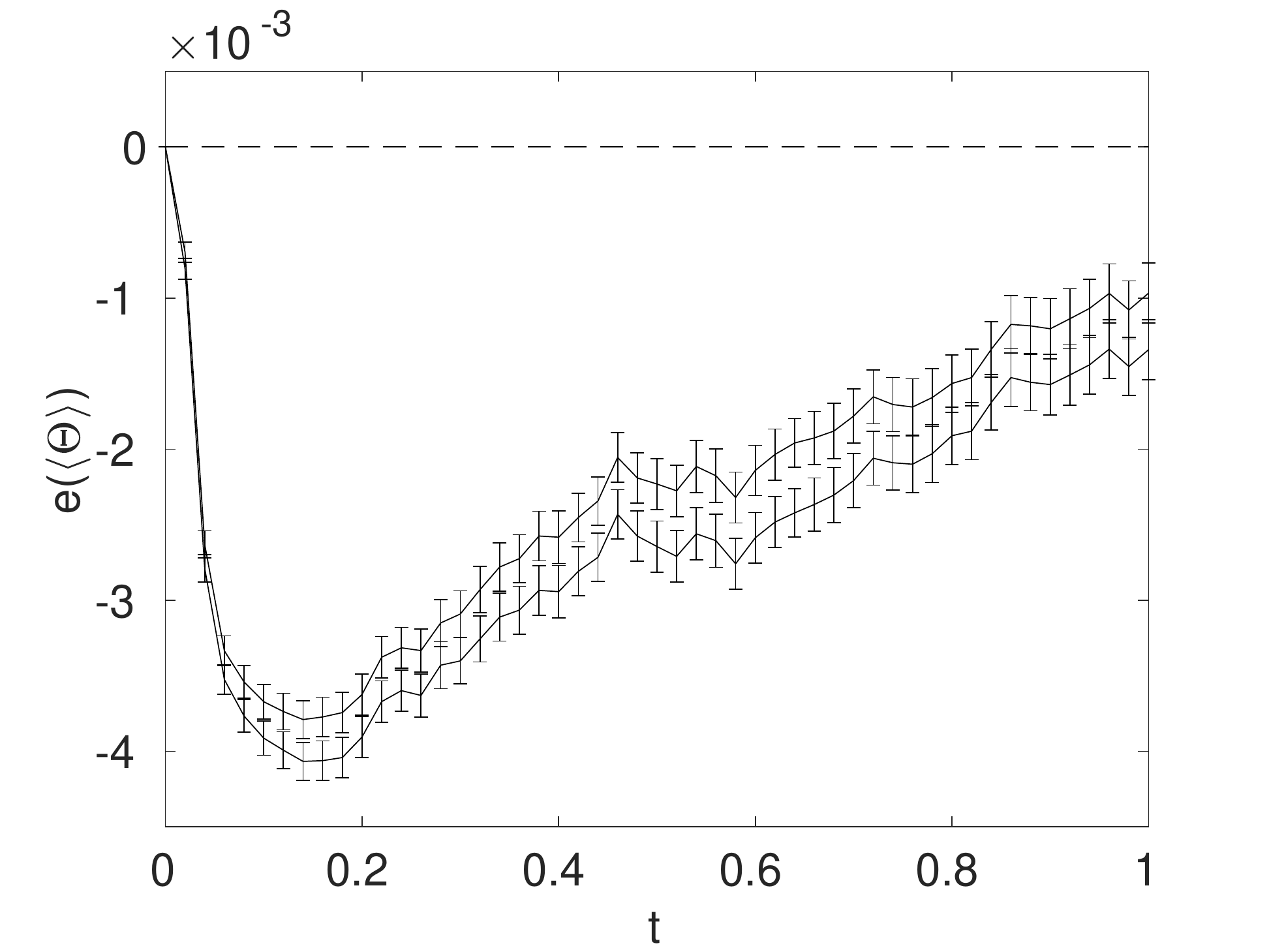}\\
\includegraphics[width=0.75\columnwidth]{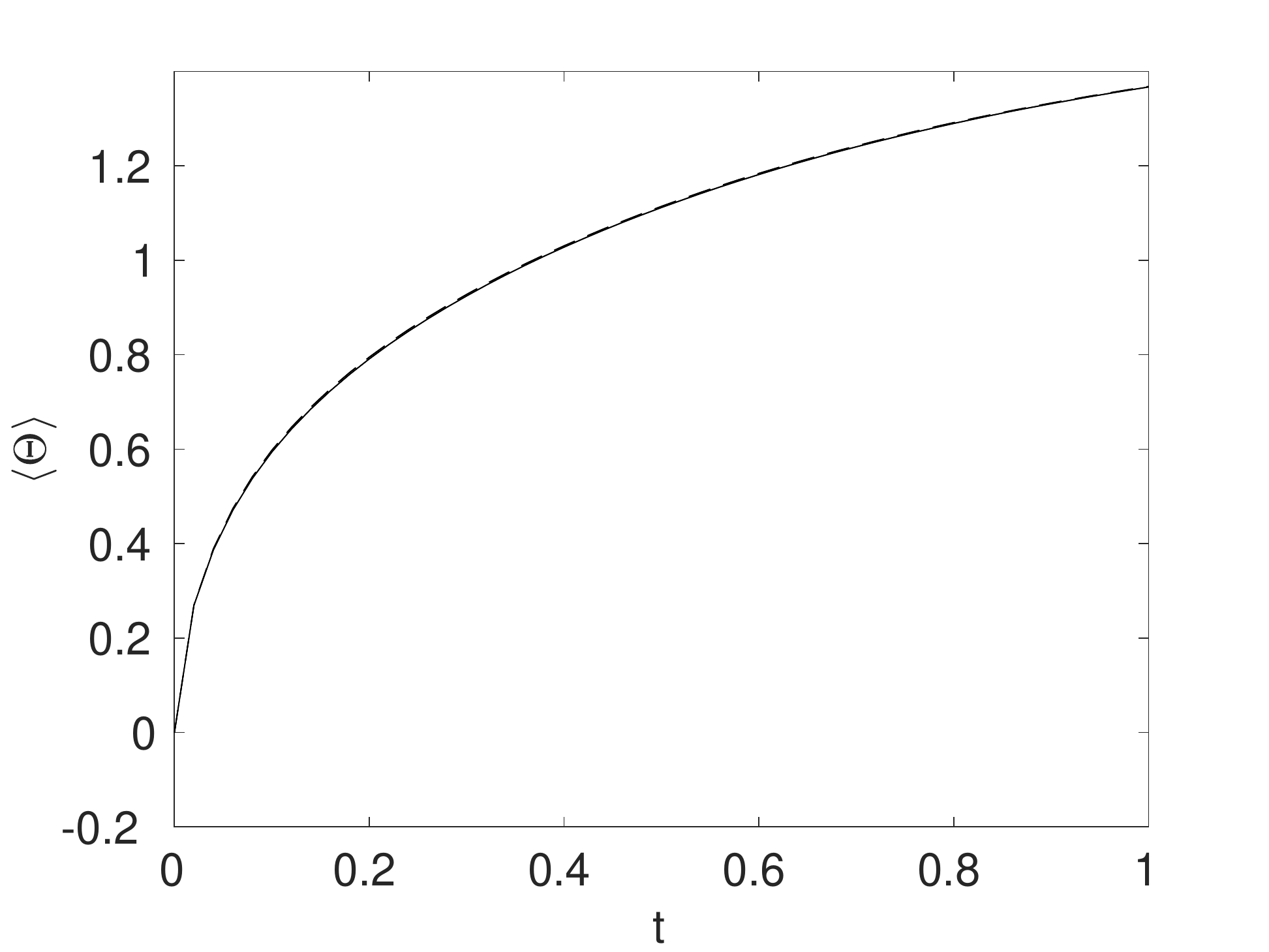}\caption{(Top) cMP mean great circle distance error, $e\left(\left\langle \Theta\right\rangle \right)$.
The two solid lines are the upper and lower sampling error bars ($\pm\sigma)$
of the projected simulation, and the error-bars indicate the sampling
errors of the intrinsic method, both with $10^{7}$ trajectories.
(Bottom) Comparison of the cMP algorithm mean great-circle distance
(solid lines) to a simulation of diffusion on the spheroid in intrinsic
coordinates (dotted lines). \label{fig:Ellipsoid-hMP1}}
\end{figure}

\begin{figure}[H]
\begin{centering}
\includegraphics[width=0.75\columnwidth]{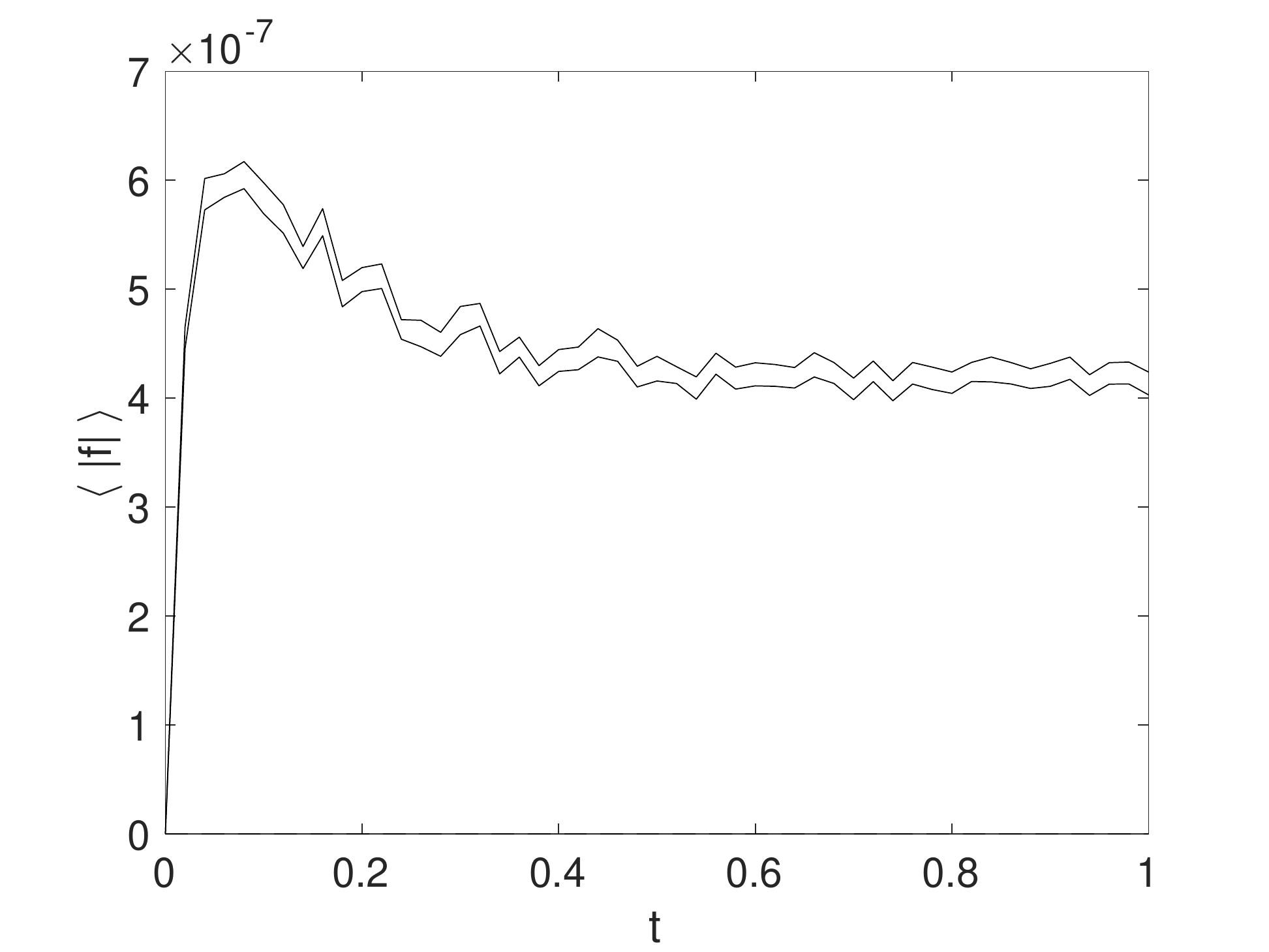}
\par\end{centering}
\caption{Constraint error of the cMP algorithm simulation of diffusion on the
spheroid. The two solid lines are the upper and lower sampling error
bars ($\pm\sigma)$ of the projected simulation. \label{fig:Ellipsoid-hMP2}}
\end{figure}

The initial point in the simulations was at $\theta=\phi=$1, which
is a large distance from any singularity. This was allowed to diffuse
following the spheroidal diffusion equations for a time-interval of
$t_{max}=1$. The results are shown in Figs (\ref{fig:Ellipsoid-hEP},
\ref{fig:Ellipsoid-MP}, \ref{fig:Ellipsoid-hMP1} \& \ref{fig:Ellipsoid-hMP2}),
where $c=0.25$ was chosen. The measure of distance used for these
figures is the great circle distance, defined as
\begin{equation}
\Theta(\bm{x},\bm{y})=\cos^{-1}\left(\bm{x}^{T}G\bm{y}\right).
\end{equation}
The results given in Table (\ref{tab:Comparison-of-error-spheroid})
show that the cMP algorithm gives excellent agreement between the
intrinsic and the projected SDE results.

\begin{table}[H]
\centering{}%
\begin{tabular}{|c|c|c|c|c|}
\hline 
Function & $\Delta t$ & cEP & tMP & \textbf{cMP}\tabularnewline
\hline 
\hline 
$\left\langle \Theta(\bm{x},\bm{x_{0}})^ {}\right\rangle $ & $0.02$ & $0.097$ & $0.1$ & $1.5\times10^{-2}$\tabularnewline
\hline 
 & $0.01$ & $6.4\times10^{-2}$ & $3.7\times10^{-2}$ & $4.2\times10^{-3}$\tabularnewline
\hline 
$\left\langle \left|f\right|\right\rangle $ & $0.02$ & $1.08\times10^{-4}$ & $0.4$ & $1.3\times10^{-5}$\tabularnewline
\hline 
 & $0.01$ & $9.9\times10^{-6}$ & $0.19$ & $6\times10^{-7}$\tabularnewline
\hline 
\end{tabular}\caption{Comparison of maximum error of combined Euler projection (cEP), tangential
midpoint projection (tMP), and combined midpoint projection (cMP)
methods for diffusion on an spheroidal surface, using time-steps of
$0.02$ and $0.01$ with $10^{7}$ parallel trajectories. Distance
sampling errors are at most $2.4\times10^{-4}$. \label{tab:Comparison-of-error-spheroid}}
\end{table}

Owing to the strong curvature of spheroid, we use smaller step-sizes
than the previous examples. For step-sizes of $0.01,$ the computational
advantage of the cMP method in this case is $16$ times reduction
in distance error compared to cEP, and $1.6\times10^{5}$ times reduction
in constraint error compared to tMP.

For all figures, the two solid lines are the upper and lower $\pm\sigma$
bounds from sampling errors in the differences, and a step size of
$0.01$ is used. 

\subsubsection{Combined Euler projection algorithm (cEP)}

Because of the final projection employed, this algorithm shows a typical
behavior of relatively high distance errors, but good ability to maintain
the trajectory on the manifold.

\subsubsection{Tangential midpoint projection algorithm (tMP)}

Because of the midpoint tangential projection employed, this algorithm
shows improved convergence with step-size, but is unable maintain
the trajectory on the manifold.

\subsubsection{Combined midpoint projection algorithm (cMP)}

Using the combined tangential projection employed, this algorithm
shows greatly improved convergence with step-size, and is able to
maintain the trajectory on the manifold.

\subsection{Hyperboloid}

The next example we consider is that of the one-sheeted hyperboloid,
defined by the constraint equation \citep{gray2017modern}
\begin{equation}
x^{2}/a^{2}+y^{2}/b^{2}-z^{2}/c^{2}-1=0.
\end{equation}
and the intrinsic coordinates:
\begin{align}
x & =\cosh v\cos\theta,\nonumber \\
y & =\cosh v\sin\theta,\nonumber \\
z & =c\sinh v.
\end{align}
We choose $a=b=1$ and $c=0.25$ to simplify the metric tensor, which
is required to obtain the intrinsic SDEs for diffusion on this manifold.
Details are given in the Appendix.

\noindent The corresponding intrinsic equations are
\begin{align}
\dot{v} & =\frac{\tanh v}{c^{2}-1+(c^{2}+1)\cosh2v}+\frac{\xi^{v}}{\sqrt{\sinh^{2}v+c^{2}\cosh^{2}v}},\nonumber \\
\dot{\theta} & =\frac{\xi^{\theta}}{\cosh v}.
\end{align}

\begin{figure}[H]
\begin{centering}
\includegraphics[width=0.75\columnwidth]{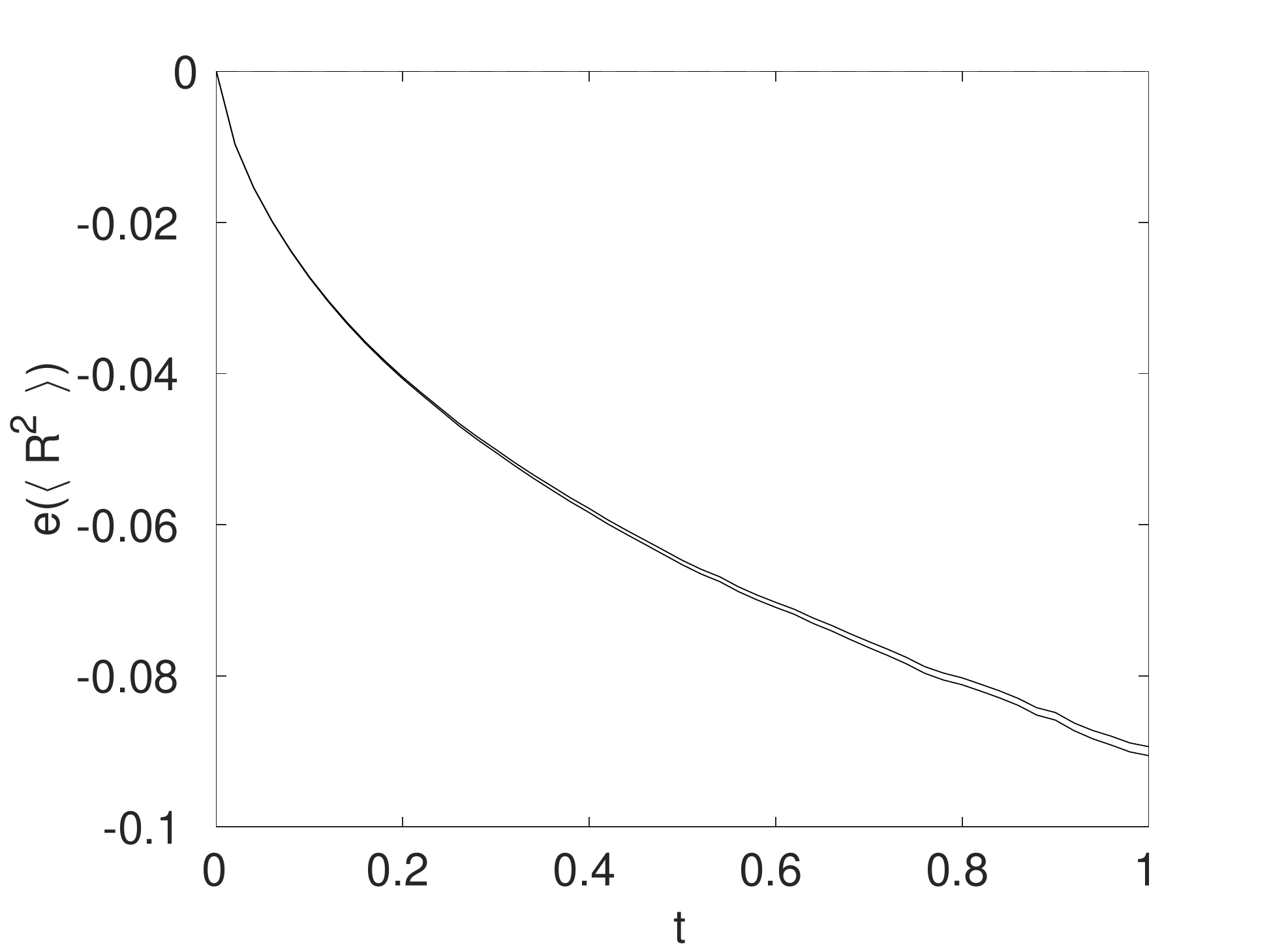}
\par\end{centering}
\centering{}\caption{(Top): cEP mean square distance error, $e\left(\left\langle R^{2}\right\rangle \right)$
is plotted. (Bottom) Constraint error of the cEP method for diffusion
on a hyperboloid. \label{fig:Hyperboloid-hEP}}
\end{figure}

\begin{figure}[H]
\begin{centering}
\includegraphics[width=0.75\columnwidth]{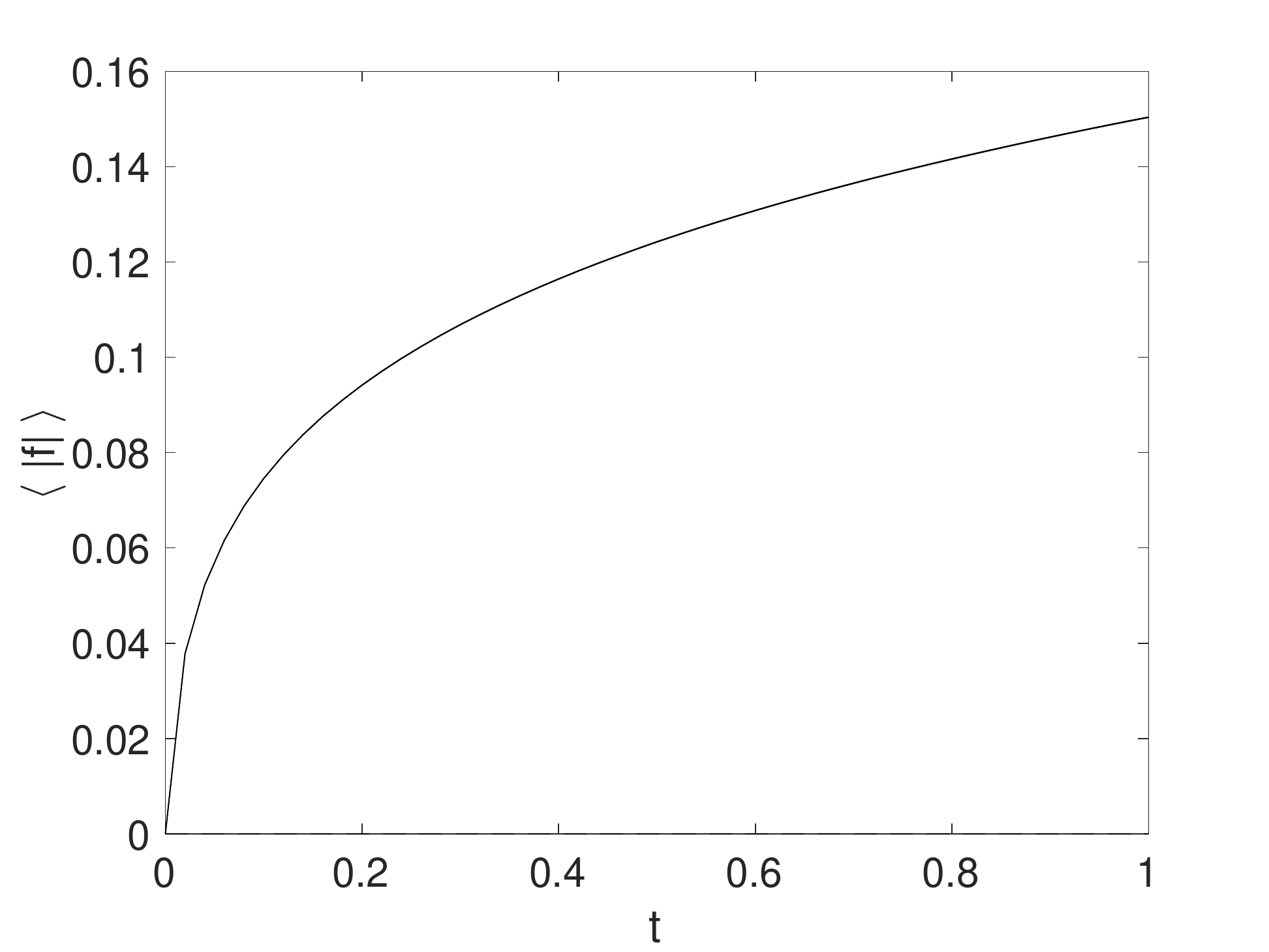}
\par\end{centering}
\caption{(Top): tMP mean square distance error, $e\left(\left\langle R^{2}\right\rangle \right)$
is plotted. (Bottom) Constraint error of the tMP algorithm for diffusion
on the hyperboloid. \label{fig:Hyperboloid-MP}}
\end{figure}

\begin{figure}[H]
\begin{centering}
\includegraphics[width=0.75\columnwidth]{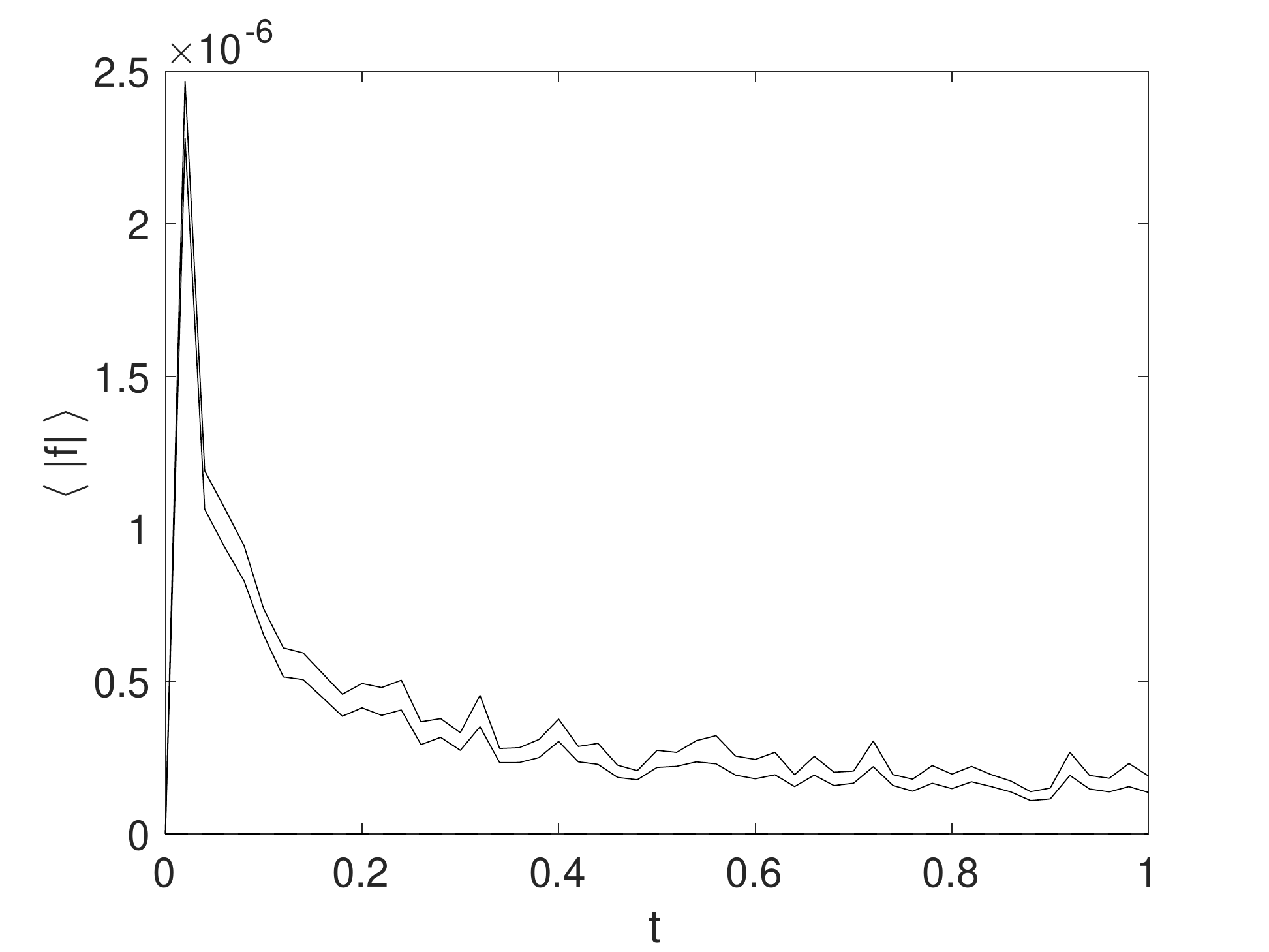}
\par\end{centering}
\caption{Constraint error of the cMP algorithm simulation of diffusion on the
hyperboloid. \label{fig:Hyperboloid-hMP2}}
\end{figure}

\begin{figure}[H]
\begin{centering}
\includegraphics[width=0.75\columnwidth]{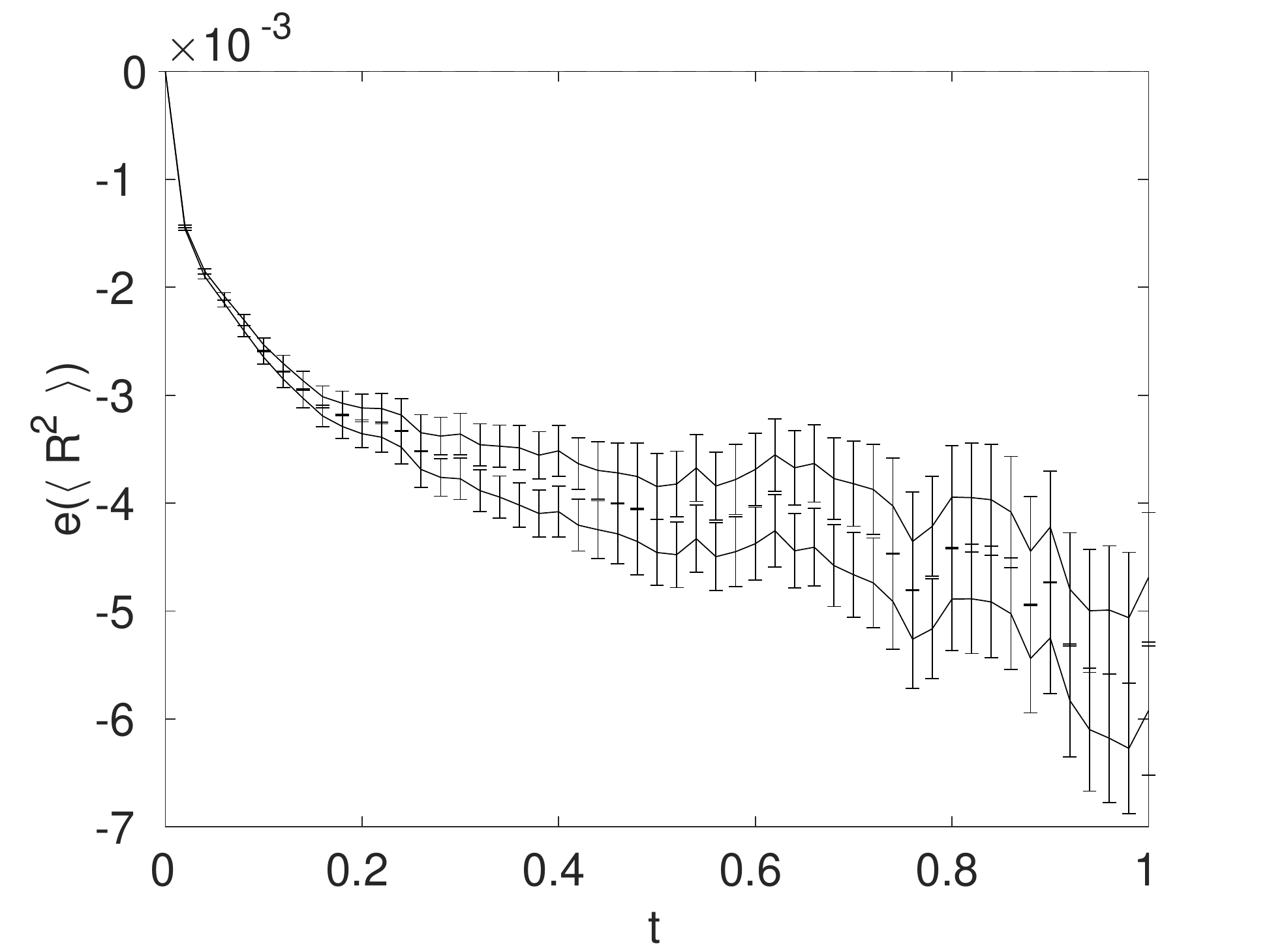}
\par\end{centering}
\centering{}\includegraphics[width=0.75\columnwidth]{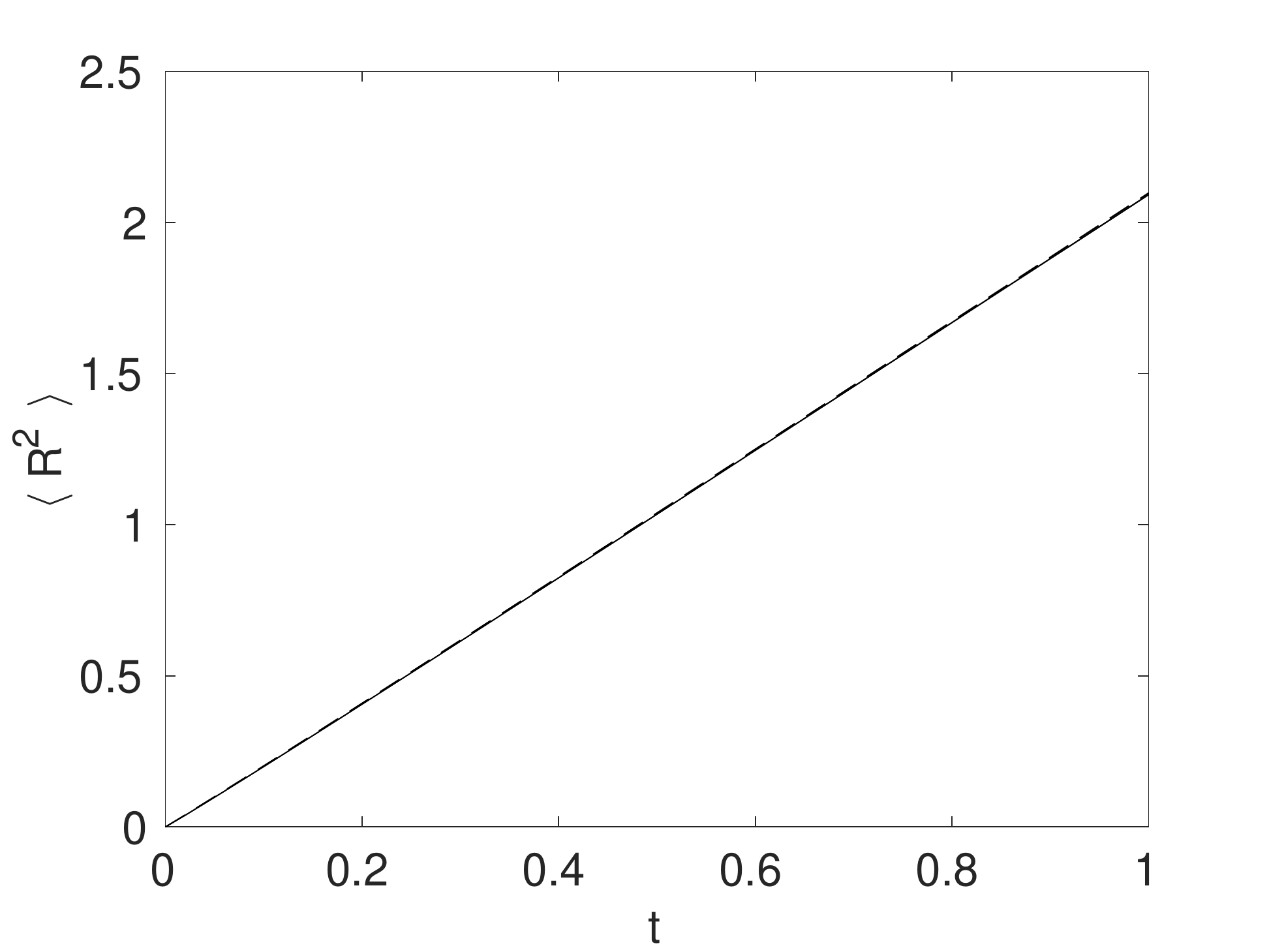}\caption{(Top) cMP mean square distance error is plotted. (Bottom) Comparison
of the cMP algorithm mean square Euclidean distance (solid lines)
to a simulation of diffusion on the hyperboloid in intrinsic coordinates
(dotted lines). \label{fig:Hyperboloid-hMP1}}
\end{figure}

Applying the three different projection algorithms for SDEs yields
Figs (\ref{fig:Hyperboloid-hEP}, \ref{fig:Hyperboloid-MP}, \ref{fig:Hyperboloid-hMP1}
\& \ref{fig:Hyperboloid-hMP2}), where we use the Euclidean distance
\citep{cartan1971differential} measure,
\begin{equation}
R(\bm{x},\bm{y})=\|\bm{x}-\bm{y}\|.
\end{equation}

The initial point in the simulations was at $\theta=v=0$ . This was
allowed to diffuse following the spheroidal diffusion equations for
a time-interval of $t_{max}=1$. Similar to the spheroid, there is
almost exact agreement between the intrinsic and hybrid midpoint projected
SDE algorithms, with poor convergence in the other two cases.

\begin{table}[H]
\centering{}%
\begin{tabular}{|c|c|c|c|c|}
\hline 
Function & $\Delta t$ & cEP & tMP & \textbf{cMP}\tabularnewline
\hline 
\hline 
$\left\langle \left|\bm{x}-\bm{x}_{0}\right|^{2}\right\rangle $ & $0.02$ & $0.14$ & $1.4\times10^{-2}$ & $1.4\times10^{-2}$\tabularnewline
\hline 
 & $0.01$ & $9\times10^{-2}$ & $1.9\times10^{-3}$ & $4.7\times10^{-3}$\tabularnewline
\hline 
$\left\langle \left|f\right|\right\rangle $ & $0.02$ & $2.9\times10^{-4}$ & $0.29$ & $3\times10^{-4}$\tabularnewline
\hline 
 & $0.01$ & $1.8\times10^{-5}$ & $0.15$ & $2.4\times10^{-6}$\tabularnewline
\hline 
\end{tabular}\caption{Comparison of maximum error of combined Euler projection (cEP), tangential
midpoint projection (tMP), and combined midpoint projection (cMP)
methods for diffusion on a hyperboloidal surface, using time-steps
of $0.02$ and $0.01$ with $10^{7}$ parallel trajectories. Distance
sampling errors are at most $6.9\times10^{-4}$. \label{tab:Comparison-of-error-hyperboloid}}
\end{table}

Because of the relatively strong curvature, we use step-sizes of $0.02$
and $0.01$. The above table shows that the computational advantage
of the cMP method. At a step-size of $.01$, there is about $16$
times reduction in mean squared distance error compared to cEP, and
a $10^{4}$ times reduction in constraint error compared to the tMP
algorithm.

\subsubsection{Combined Euler projection algorithm (cEP)}

As in the spheroidal case, this algorithm shows relatively poor convergence
with step-size, but is able to maintain the trajectory on the manifold,
due to the use of a normal projection.

\subsubsection{Tangential midpoint projection algorithm (tMP)}

This algorithm shows improved convergence with step-size, but is less
able to maintain the trajectory on the manifold, with very large constraint
errors.

\subsubsection{Combined midpoint projection algorithm (cMP)}

This algorithm repeats the trend in the spheroidal case, with both
improved convergence with step-size, and much greater ability to maintain
the trajectory on the manifold.

We note the almost linear growth in the mean square diffusion distance.
This is caused by our starting point at $z=c\sinh v=0$. Apart from
a scaling of the $z$-distance, in the small-$z$ regime the hyperboloidal
geometry is similar to a catenoid with $z=v$, which has an exactly
linear growth.

\subsection{Polynomial surface}

A more general $N-th$ order polynomial manifold in $n$ dimensions
was investigated, with a constraint of:
\begin{equation}
\sum_{j=1}^{n}\left[x^{j}\right]^{N}=1.\label{eq:pol}
\end{equation}
This example is applicable to arbitrary dimensions. In the three-dimensional
case, for powers $N\gg2$, this gives a nearly cubic surface with
rounded edges, corresponding to a rapidly varying curvature, as in
some biological cell walls. To give a model of transport, we include
an axial force as a drift term. 

The projected SDE can be written as in Eq (\ref{eq:noise-projection})
as:
\begin{equation}
\dot{\bm{x}}=\mathcal{P}_{\bm{x}}^{\parallel}\left(\bm{a}+\bm{\xi}\right),
\end{equation}
where $\bm{a}=\left[0,0,2z\right]$, thus modeling an outward axial
force in the $z$ direction.

Here we use a quartic constraint: $N=4$, to give an inhomogeneous,
non-spherical surface. As there are neither exact solutions nor intrinsic
variables, the reference solutions were obtained using the cMP algorithm
with much smaller step-sizes. An error comparison is shown in Table
(\ref{tab:Comparison-of-error-n-sphere-2}). 
\begin{table}[H]
\centering{}%
\begin{tabular}{|c|c|c|c|c|}
\hline 
Function & $\Delta t$ & cEP & tMP & \textbf{cMP}\tabularnewline
\hline 
\hline 
$\left\langle R^{2}\right\rangle $ & $0.1$ & $0.14$ & $0.12$ & $1.7\times10^{-2}$\tabularnewline
\hline 
 & $0.05$ & $0.08$ & $0.07$ & $9.9\times10^{-3}$\tabularnewline
\hline 
$\left\langle \left|f\right|\right\rangle $ & $0.1$ & $2.8\times10^{-4}$ & $0.3$ & $1.7\times10^{-5}$\tabularnewline
\hline 
 & $0.05$ & $1.3\times10^{-5}$ & $0.16$ & $4\times10^{-9}$\tabularnewline
\hline 
\end{tabular}\caption{Comparison of maximum error of combined Euler projection (cEP), tangential
midpoint projection (tMP), and combined midpoint projection (cMP)
methods for diffusion on a 3-dimensional, fourth-order polynomial
surface, using time-steps of $0.1$ and 0.05, with $10^{7}$ parallel
trajectories, and a duration of $t_{max}=5$. Distance sampling errors
are less than $5\times10^{-4}$. \label{tab:Comparison-of-error-n-sphere-2}}
\end{table}

As in the earlier examples, the distance errors with the cMP algorithm
are smaller by nearly an order of magnitude than with the other algorithms.

\subsection{Hypersphere}

We finally consider much higher dimensions, using an $n$-dimensional
hyper-spherical manifold with $\left|\bm{x}\right|^{2}=1$ and isotropic
unit diffusion in the original Euclidean space. If $\mathcal{P}_{\bm{x}}^{\parallel}$
is a hyperspherical surface projection, the projected SDE can be written
as in Eq (\ref{eq:noise-projection}). This has exact solutions for
the expected diffusion distance, $R=\left|\bm{x}-\bm{x}_{0}\right|$
\citep{castro2014intrinsic}:

\begin{equation}
\left\langle R^{2}\right\rangle =2\left(1-\exp\left(-\frac{n-1}{2}t\right)\right).
\end{equation}

Numerical results for errors in ten space dimensions is shown in Table
(\ref{tab:Comparison-of-error-n-sphere}). 
\begin{table}[H]
\centering{}%
\begin{tabular}{|c|c|c|c|c|}
\hline 
Function & $\Delta t$ & cEP & tMP & cMP\tabularnewline
\hline 
\hline 
$\left\langle R^{2}\right\rangle $ & $0.1$ & $0.29$ & $0.27$ & \textbf{$2.1\times10^{-2}$}\tabularnewline
\hline 
 & $0.05$ & $0.17$ & $0.11$ & \textbf{$9.8\times10^{-3}$}\tabularnewline
\hline 
$\left\langle \left|f\right|\right\rangle $ & $0.1$ & $1.0\times10^{-5}$ & $0.3$ & \textbf{$7.8\times10^{-11}$}\tabularnewline
\hline 
 & $0.05$ & $2.8\times10^{-7}$ & $0.11$ & 2.3\textbf{$\times10^{-16}$}\tabularnewline
\hline 
\end{tabular}\caption{Comparison of maximum error of combined Euler projection (cEP), tangential
midpoint projection (tMP), and combined midpoint projection (cMP)
methods for diffusion on a 10-dimensional hypersphere, using time-steps
of $0.1$ and 0.05. \label{tab:Comparison-of-error-n-sphere}}
\end{table}

\begin{figure}[H]
\centering{}\includegraphics[width=0.75\columnwidth]{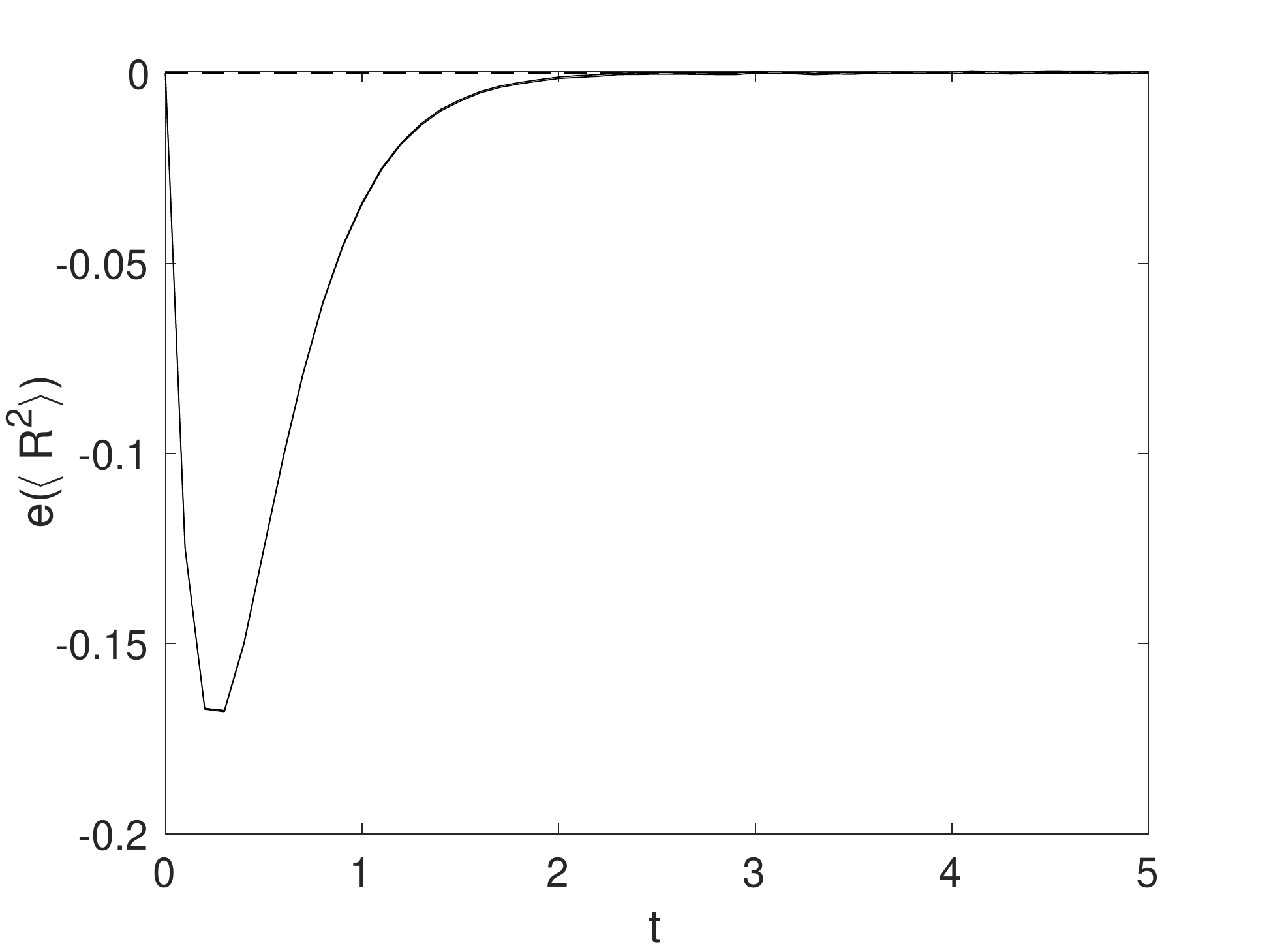}\\
\includegraphics[width=0.75\columnwidth]{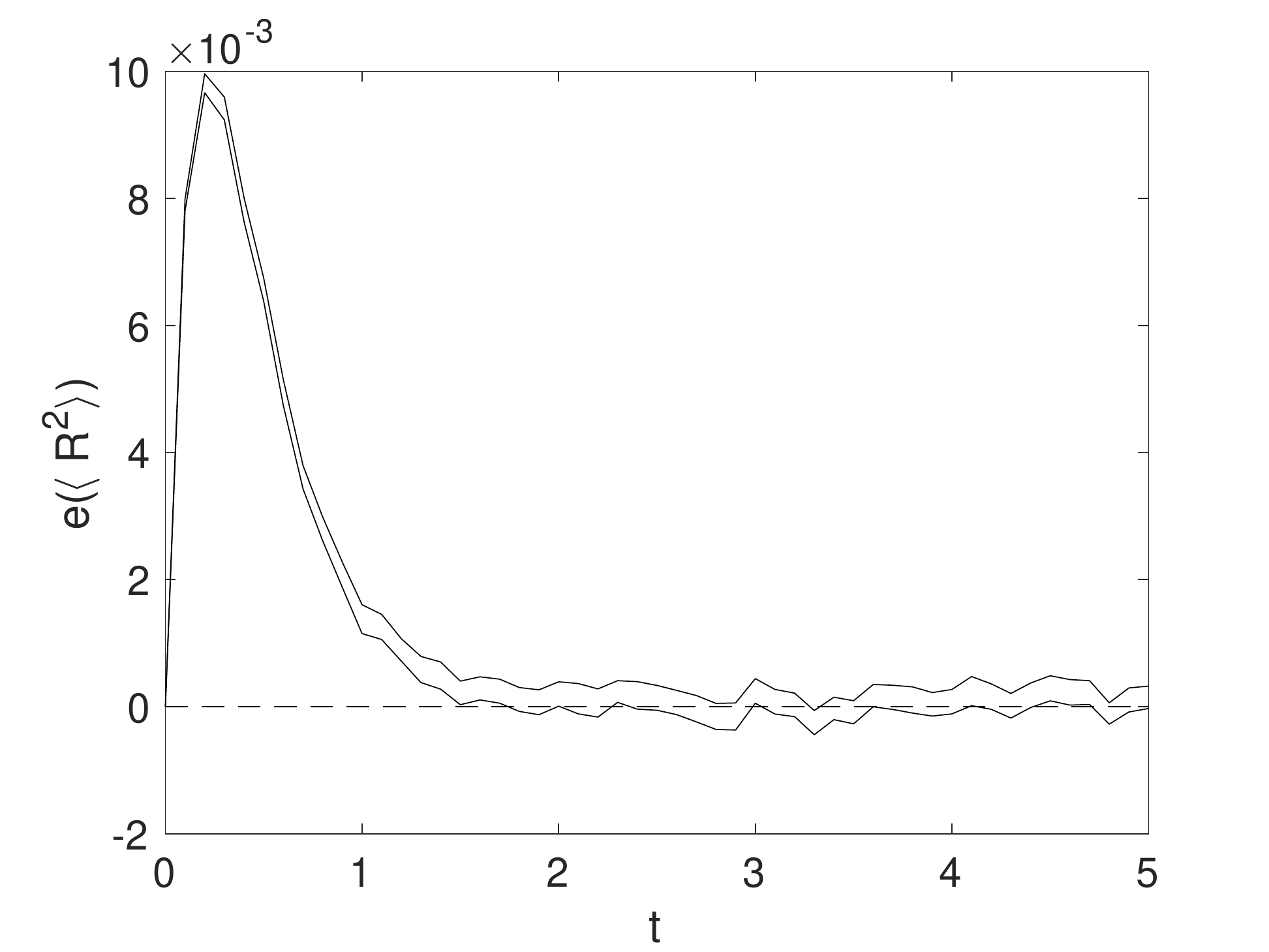}\caption{Comparing errors for the 10-sphere diffusion case with a step-size
of $0.05$. The solid lines are upper and lower sampling error bars
of the simulations. Top: combined Euler projection method (cEP). Bottom:
combined midpoint projection method (cMP). \label{fig:hypersphere}}
\end{figure}

The computational advantage of the cMP method increases in higher
dimensions. In $n=10$ dimensions, for a step-size of $0.05$, the
improvement for cMP is a factor of $17$ in mean square distance error
reductions compared to cEP, and over $10^{9}$ in constraint error
reductions compared to either alternative. In these tests, the duration
is $t_{max}=5$, with $10^{7}$ parallel trajectories. Sampling errors
are of order $2\times10^{-4}$. Graphical comparisons are shown in
Fig (\ref{fig:hypersphere}).

The ratio of errors obtained in the Kubo oscillator case, with an
8:1 improvement in error for the midpoint versus the Euler method,
is even greater in these higher-dimensional results. After correcting
for the increased diffusion in ten dimensions, these errors are within
the global error bounds of Eqs (\ref{eq:global-error-Euler}) and
(\ref{eq:Global-error-MP}).

\section{Summary of results\label{sec:Summary}}

In summary, we have shown how projected SDEs are derived using adiabatic
elimination with a constraint potential, and that the projected equation
is a Stratonovich type. We have obtained a combined midpoint projection
algorithm for projecting a stochastic equation onto a general manifold,
and compared it to earlier proposals, using numerical studies on multiple
surfaces with positive and negative curvature, having different dimensions
and constraint equations. 

The hybrid midpoint algorithm gives a distinctly improved error performance
compared to both the earlier algorithms used for comparison. The improvement
is due to a combination of tangential and normal projections, together
with a more accurate midpoint algorithm. The combined Euler projection
method cannot accurately track the changes in the projected drift
and diffusion terms during a step, as it relies on an initial estimate
of the diffusion, which changes in space due to the projection. The
midpoint projection method has an improved treatment of diffusion,
but is unable to control global error growth when trajectories move
off the manifold.

Adaptive step-size \citep{lamba2003adaptive} or higher order methods
are possible \citep{burrage1996high,burrage2000order,Kloeden1992}.
These are also suggested for projected equations \citep{abdulle2014high,laurent2021order},
and provide other alternatives.

Yet there are multiple factors \citep{opanchuk2016parallel,kiesewetter2017algorithms}
causing numerical errors. The combined midpoint projected method combines
good discretization error performance with low complexity, straightforward
parallel implementation, and reasonable speed. This is a great advantage
when it is important to reduce both sampling error and discretization
errors.

In summary, the combined midpoint projection method, which allows
for changes in the diffusion matrix and maintains the constraint,
has much greater accuracy than either method used for comparisons.
This comparative performance is proved as an error bound for the simple
case of the Kubo oscillator, and demonstrated numerically in more
complex examples. This implies that, for identical error performance,
the cMP method can use larger steps, and is up to an order of magnitude
more efficient than the cEP algorithm. 

\section*{Acknowledgements}

This work was funded through an Australian Research Council Discovery
Project Grant DP190101480, and a grant from NTT Research. PDD acknowledges
the hospitality of the Aspen Center for Physics, supported by NSF
grant PHY-1607611, the Institute for Atomic and Molecular Physics
(ITAMP) at Harvard University, supported by the NSF, and the Joint
Institute for Laboratory Astrophysics at University of Colorado.

\section*{Appendix: Spheroidal and hyperboloidal intrinsic stochastic equations}

\subsection*{Spheroidal diffusion}

As an example in the main text, a 2-spheroid may be embedded in $\mathbb{R}^{3}$
by the constraint equation \citep{gray2017modern}
\begin{equation}
x^{2}+y^{2}+z^{2}/c^{2}=1.
\end{equation}
Applying Eq. (\ref{eq:pullback}), one obtains the inverse metric
tensor 
\begin{equation}
(g^{\mu\nu})=\text{\ensuremath{\mathrm{diag}\left(\frac{1}{\cos^{2}\theta+c^{2}\sin^{2}\theta},\frac{1}{\sin^{2}\theta}\right)},}\label{eq:Diff tensor}
\end{equation}
where the intrinsic coordinates $\theta$ and $\phi$ are related
to the extrinsic coordinates as \citep{gray2017modern}
\begin{align}
x & =\sin\theta\cos\phi,\nonumber \\
y & =\sin\theta\sin\phi,\nonumber \\
z & =c\cos\theta,
\end{align}
with $0\leq\theta\leq\pi$ and $0\leq\phi<2\pi.$ 

From the diffusion tensor in Eq. (\ref{eq:Diff tensor}) and the relation
in Eq. (\ref{eq:SelfAdjointDrift}), one can find the self adjoint
drift as,

\begin{align}
\alpha_{a}^{\theta} & =\frac{\mathrm{cos}\theta\left[\mathrm{sin}\theta\left(\mathrm{c^{2}}-1\right)+\frac{1}{2\mathrm{sin}\theta}\right]}{\left(\mathrm{cos^{2}}\theta+\mathrm{c^{2}}\sin^{2}\theta\right)^{2}},\\
\alpha_{a}^{\phi} & =0.
\end{align}
This can be converted to a Stratonovich drift by making use of Eq.
(\ref{eq:StratDrift}):

\begin{align}
\alpha_{\text{}}^{\theta} & =\frac{\mathrm{cot}\theta}{2\left(\mathrm{cos^{2}}\theta+\mathrm{c^{2}}\sin^{2}\theta\right)},\\
\alpha_{\text{}}^{\phi} & =0.
\end{align}

This in turn can be utilized to find the intrinsic Stratonovich SDE
as:

\begin{align}
\dot{\theta} & =-\frac{\cot\theta}{(c^{2}-1)\cos2\theta-(1+c^{2})}+\frac{\xi^{\theta}}{\sqrt{\cos^{2}\theta+c^{2}\sin^{2}\theta}}\nonumber \\
\dot{\phi} & =\frac{\xi^{\phi}}{\sin\theta}.
\end{align}

\subsection*{Hyperboloidal diffusion}

\noindent As another example used in the main text, the one-sheeted
hyperboloid is defined by the constraint \citep{gray2017modern,berger2012differential}
\begin{equation}
x^{2}+y^{2}-z^{2}/c^{2}=1.
\end{equation}
A convenient set of intrinsic coordinates is given by \citep{gray2017modern}
\begin{align}
x & =\cosh v\cos\theta,\nonumber \\
y & =\cosh v\sin\theta,\nonumber \\
z & =c\sinh v,
\end{align}
where $-\infty<v<\infty$ and $0\leq\theta<\pi.$ This yields the
inverse metric
\begin{equation}
(g^{\mu\nu})=\mathrm{diag}\left(\frac{1}{\sinh^{2}v+c^{2}\cosh^{2}v},\frac{1}{\cosh^{2}v}\right),
\end{equation}
and making using of the above equation, (\ref{eq:SelfAdjointDrift})
and (\ref{eq:StratDrift}) one obtains an intrinsic Stratonovich SDE:
\begin{align}
\dot{v} & =\frac{\tanh v}{c^{2}-1+(c^{2}+1)\cosh2v}+\frac{\xi^{v}}{\sqrt{\sinh^{2}v+c^{2}\cosh^{2}v}},\nonumber \\
\dot{\theta} & =\frac{\xi^{\theta}}{\cosh v}.
\end{align}

\bibliographystyle{apsrev4-1}

\end{document}